\documentclass{article}
\usepackage{graphicx}
\usepackage{amsmath}
\usepackage{amsfonts}
\usepackage{amssymb}
\newtheorem{theorem}{Theorem}

\newtheorem{corollary}[theorem]{Corollary}

\newtheorem{definition}[theorem]{Definition}

\newtheorem{lemma}[theorem]{Lemma}

\newtheorem{proposition}[theorem]{Proposition}

\begin{document}

\title{Riemannian Geometry \\over different\\Normed Division Algebras}
\author{Naichung Conan Leung}
\maketitle
\tableofcontents

\label{1Sec Intro}\pagebreak 

\textbf{Abstract}: We develop a unifed theory to study geometry of manifolds
with different holonomy groups. They are classified by (1) real, complex,
quaternion or octonion number they are defined over and (2) being special or
not. Specialty is an orientation with respect to the corresponding normed
algebra $\mathbb{A}$. For example, special Riemannian $\mathbb{A}$-manifolds
are oriented Riemannian, Calabi-Yau, Hyperk\"{a}hler and $G_{2}$-manifolds respectively.

For vector bundles over such manifolds, we introduce (special) $\mathbb{A}%
$-connections. They include holomorphic, Hermitian Yang-Mills, Anti-Self-Dual
and Donaldson-Thomas connections. Similarly we introduce (special) $\frac
{1}{2}\mathbb{A}$-Lagrangian submanifolds as maximally real submanifolds. They
include (special) Lagrangian, complex Lagrangian, Cayley and (co-)associative submanifolds.

We also discuss geometric dualities from this viewpoint: Fourier
transformations on $\mathbb{A}$-geometry for flat tori and a conjectural SYZ
mirror transformation from (special) $\mathbb{A}$-geometry to (special)
$\frac{1}{2}\mathbb{A}$-Lagrangian geometry on mirror special $\mathbb{A}$-manifolds.

\pagebreak 

\section{Introduction}

It is well-known that a Riemannian metric $g$ on a manifold $M$ determines a
unique torsion free Riemannian connection on its tangent bundle, called the
Levi-Civita connection. For a generic metric $g$, its holonomy group
$Hol\left(  g\right)  $ equals $O\left(  m\right)  $ with $m=\dim M$. The size
of the holonomy group is inversely proportional to the amount of geometric
structures $M$ possesses. For example $Hol\left(  g\right)  \subset U\left(
n\right)  $, with $m=2n$, is equivalent to $M$ being a K\"{a}hler manifold.
When we further restrict the holonomy group to $SU\left(  n\right)  $, we
obtain a Calabi-Yau manifold and they are the central objects of interest in
mirror symmetry. Recently, M-theory suggests that the geometry of seven
dimensional manifolds with $Hol\left(  g\right)  \subset G_{2}$ has even
richer geometry. Other holonomy groups, like $Sp\left(  n\right)  $ for
hyperk\"{a}hler manifolds, are also very important in modern geometry. A
complete classification of all possible holonomy groups has been obtained by
Berger \cite{Berger} many years ago.

In this paper we are going to study all these geometries from a unified point
of view. Namely we analyze geometries as they are defined over $\mathbb{R}$,
$\mathbb{C}$, $\mathbb{H}$ or $\mathbb{O}$, the four normed division algebras
$\mathbb{A}$. In a sense this approach is very natural, as \textit{metric
geometry }should be defined over \textit{metric algebras}.

A Riemannian manifold is called an\textit{ }$\mathbb{A}$\textit{-manifold }if
its holonomy group is inside $G_{\mathbb{A}}\left(  n\right)  $, the group of
twisted isomorphisms of $\mathbb{A}^{n}$. We can identify these manifolds for
various $\mathbb{A}$ as Riemannian manifolds, K\"{a}hler manifolds,
quaternionic K\"{a}hler manifolds and $Spin\left(  7\right)  $-manifolds
respectively. There is also a notion of an $\mathbb{A}$-orientation and it
defines a subgroup $H_{\mathbb{A}}\left(  n\right)  $ in $G_{\mathbb{A}%
}\left(  n\right)  $ consisting of special twisted isomorphisms of
$\mathbb{A}^{n}$. Corresponding manifolds are called \textit{special
}$\mathbb{A}$\textit{-manifolds,} and they are oriented Riemannian manifolds,
Calabi-Yau manifolds, hyperk\"{a}hler manifolds and $G_{2}$-manifolds
respectively. Notice that these are precisely all possible holonomy groups for
a Riemannian manifold which is not locally symmetric. In this classification,
manifolds with $G_{2}$ holonomy group are those with richest geometric structures.

Unlike real or complex manifolds, $\mathbb{H}$- and $\mathbb{O}$-manifolds do
not have many functions. Their geometries are reflected by submanifolds and
bundles over them. To define them in a unified way, we note that holonomy
groups $G_{\mathbb{A}}\left(  n\right)  $ and $H_{\mathbb{A}}\left(  n\right)
$ define natural subbundles $\mathbf{g}_{\mathbb{A}}\left(  T_{M}\right)  $
and $\mathbf{h}_{\mathbb{A}}\left(  T_{M}\right)  $ in $\Lambda^{2}T_{M}%
^{\ast}$. For example when $\mathbb{A}=\mathbb{C}$, i.e. $M$ is a K\"{a}hler
manifold, we have $\mathbf{g}_{\mathbb{C}}\left(  T_{M}\right)  =\Lambda
^{1,1}\left(  M\right)  $ and $\mathbf{h}_{\mathbb{C}}\left(  T_{M}\right)
=\Lambda_{0}^{1,1}\left(  M\right)  $.

A connection $D_{E}$ on a vector bundle $E$ over $M$ is called an $\mathbb{A}%
$-\textit{connection} (resp. \textit{special }$\mathbb{A}$\textit{-connection}%
) if its curvature tensor $F_{E}$ lies inside $\mathbf{g}_{\mathbb{A}}\left(
T_{M}\right)  \otimes ad\left(  E\right)  $ (resp. $\mathbf{h}_{\mathbb{A}%
}\left(  T_{M}\right)  \otimes ad\left(  E\right)  $). In the K\"{a}hler case,
such a connection $D_{E}$ is a holomorphic connection, i.e. $F_{E}^{0,2}%
=F_{E}^{2,0}=0$ (resp. Hermitian Yang-Mills connection). Special $\mathbb{A}%
$-connections, with $\mathbb{A}\neq\mathbb{R}$, are always absolute minimum
for the Yang-Mills energy functional, as it will be explained in terms of
bundle calibrations.

A natural class of submanifolds in any $\mathbb{A}$-manifold consists of
$\mathbb{A}$-submanifolds. However, there is another natural class of
submanifolds in the middle dimension, which plays the role of
\textit{decomplexifying} $M$, and they are called $\frac{1}{2}\mathbb{A}%
$-\textit{Lagrangian submanifolds}. For example they include Lagrangian
submanifolds in K\"{a}hler manifolds and Cayley submanifolds in $Spin\left(
7\right)  $-manifolds. Using the $\mathbb{A}$-orientation on a special
$\mathbb{A}$-manifold, we also have the notion of \textit{special} $\frac
{1}{2}\mathbb{A}$-\textit{Lagrangian submanifolds of type I }or \textit{type
II}. They can be identified as special Lagrangian submanifolds (with phase
angle $0$ or $\pi/2$) in Calabi-Yau manifolds, complex Lagrangian submanifolds
in hyperk\"{a}hler manifolds, associative submanifolds and coassociative
submanifolds in $G_{2}$-manifolds. As in the bundle case, special $\frac{1}%
{2}\mathbb{A}$-Lagrangian submanifolds are absolute minimum for the volume
functional as they are all volume calibrated.

A good notion of a \textit{global} \textit{decomplexification} of $M$ is a
fibration with a section on $M$ by (special) $\frac{1}{2}\mathbb{A}%
$-Lagrangian submanifolds, possibly singular. For example, in the theory of
geometric quantization of symplectic manifolds, a real polarization is merely
a smooth Lagrangian fibration with a section.

To define and study the geometry of any (special) $\mathbb{A}$-manifold $M$,
we need to couple submanifolds $C$ in $M$ with connections $D_{E}$ over $C$
and we call any such pair $\left(  C,D_{E}\right)  $ a \textit{cycle}. We have
(i) a \textit{(special)} $\mathbb{A}$\textit{-cycle} consists of an
$\mathbb{A}$-submanifold and a (special) $\mathbb{A}$-connection over it and
(ii) a \textit{(special)} $\frac{1}{2}\mathbb{A}$\textit{-Lagrangian cycle}
consists of a (special) $\frac{1}{2}\mathbb{A}$-Lagrangian submanifold and a
special $\mathbb{A}$-connection over it. For instance, in the K\"{a}hler case,
a $\mathbb{C}$-cycle is a holomorphic bundle over a complex submanifold in
$M$, in particular it is a coherent sheaf on $M$. Such a $\mathbb{C}$-cycle is
special if the bundle carries a Hermitian Yang-Mills connection. Similarly, a
(special) $\mathbb{R}$-Lagrangian cycle is a flat bundle over a (special)
Lagrangian submanifold in the Calabi-Yau manifold $M$. The mirror symmetry
conjecture says that these two types of geometries can be transformed to each
other on mirror Calabi-Yau manifolds.

Such duality transformations play very important roles both in mathematical
physics and geometry. A basic ingredient is the Fourier transformation. We
will recall how it transforms the $\mathbb{A}$-geometry on a flat torus over
$\mathbb{A}$ to the $\mathbb{A}$-geometry on its dual torus. For a special
$\mathbb{A}$-manifold $M$ with a special $\frac{1}{2}\mathbb{A}$-Lagrangian
fibration and a section, we also discuss briefly the SYZ mirror duality which
transforms the (special) $\mathbb{A}$-geometry of $M$ to the (special)
$\frac{1}{2}\mathbb{A}$-Lagrangian geometry of its mirror manifold. This can
be viewed as a fiberwise Fourier transformation along $\frac{1}{2}\mathbb{A}%
$-Lagrangian fibrations (see \cite{SYZ}, \cite{Gross TopMS}, \cite{LYZ},
\cite{KontSoi}, \cite{Le Geom MS} in the Calabi-Yau case, \cite{Le Fourier HK}
in the hyperk\"{a}hler case and \cite{Acharya}, \cite{GYZ}, \cite{LL} in the
$G_{2}$-manifolds case).

In the last section, we give several remarks and questions on related aspects
of geometries over different normed division algebras $\mathbb{A}$.

\section{\label{1Sec Sp A mfd}(Special) Riemannian $\mathbb{A}$-manifolds}

In this section we define Riemannian manifolds over different normed division
algebras. As we will see, \textit{all} possible holonomy groups are naturally
aroused from manifolds defined over $\mathbb{R}$, $\mathbb{C}$, $\mathbb{H}$
or $\mathbb{O}$; unless they are locally symmetric spaces. This gives us a
unified way to look at all the seemingly unrelated branches of geometry in
mathematics. We first recall some basic facts about normed division algebras
$\mathbb{A}$ (see e.g. \cite{HL}), then we will define the key notion: twisted isomorphisms.

\subsection{\label{1Sec Norm alg}Normed division algebra, $\mathbb{A}$}

\begin{definition}
A normed algebra $\mathbb{A}$ is a finite dimensional real algebra with a unit
$1$ and a norm $\left\|  \cdot\right\|  $ satisfying $\left\|  a\cdot
b\right\|  =$ $\left\|  a\right\|  $ $\left\|  b\right\|  $ for any
$a,b\in\mathbb{A}$.
\end{definition}

There are exactly four of them, namely the real $\mathbb{R}$, the complex
$\mathbb{C}$, the quaternion $\mathbb{H}$ and the octonion (or Cayley)
$\mathbb{O}$ numbers. Each can be interpreted as the complexification of the
previous one, the so-called Cayley-Dickson process: suppose $\mathbb{A}$ is
any algebra with conjugation $\ast$,\footnote{We also denote $\ast a$ as
$\bar{a}$.} we define an algebra structure on $\mathbb{B}=\mathbb{A}%
\oplus\mathbb{A}$ as follows,
\begin{align*}
\left(  a,b\right)  \left(  c,d\right)   &  =\left(  ac-db^{\ast},a^{\ast
}d+cb\right) \\
\left(  a,b\right)  ^{\ast}  &  =\left(  a^{\ast},-b\right)  \text{.}%
\end{align*}
This process will construct $\mathbb{C}$ from $\mathbb{R}$ (and we write
$\mathbb{R}=\frac{1}{2}\mathbb{C}$) and so on.

The following properties of a normed algebra will be needed in this article:
(i) $\left\langle xy,z\right\rangle =\left\langle x,z\bar{y}\right\rangle $,
(ii) $\left\langle xy,zy\right\rangle =\left\langle x,z\right\rangle \left|
y\right|  ^{2}$ and (iii) $\left(  xy\right)  y=x\left(  y^{2}\right)  $ for
any $x,y,z\in\mathbb{A}$. Furthermore, $\mathbb{A}$ is always a division algebra.

Each time we complexify $\mathbb{A}$, we loss some nice properties: (1)
$\mathbb{C}$ is not real, (2) $\mathbb{H}$ is not commutative, (3)
$\mathbb{O}$ is not associative and lastly $\mathbb{O\oplus O}$ is no longer
normed. As a result, $\mathbb{H}^{n}$ is only a bi-module of $\mathbb{H}$ and
not a vector space because $\mathbb{H}$ is not a field. Furthermore,
$\mathbb{O}$ does not even act on $\mathbb{O}^{n}$, with $n\geq2$, because of
the non-associativity. We call $\mathbb{A}^{n}$ a \textit{linear }$\mathbb{A}%
$-\textit{space}\textbf{\ }of rank $n$, and its real dimension is $m=2^{a}n$
with $a=0,1,2$ or $3$. An inner product on $\mathbb{A}^{n}$ always refers to
one satisfying
\[
\left\langle u\cdot x,v\cdot x\right\rangle =\left\langle u,v\right\rangle
\left|  x\right|  ^{2},
\]
for any $u,v\in V$ and $x\in\mathbb{A}$.

\bigskip

To define $\mathbb{A}$-manifolds, the nonlinear analog of linear $\mathbb{A}%
$-spaces, we first need to define twisted isomorphisms of $\mathbb{A}^{n}$.
Recall on a normed linear $\mathbb{A}$-space $V\cong\mathbb{A}^{n}$, an
automorphism of $V$ is a real linear isomorphism $\phi:V\rightarrow V$
satisfying
\begin{align*}
\left\langle \phi\left(  u\right)  ,\phi\left(  v\right)  \right\rangle  &
=\left\langle u,v\right\rangle ,\\
\phi\left(  vx\right)   &  =\phi\left(  v\right)  x\text{,}%
\end{align*}
for any $u,v\in V$ and $x\in\mathbb{A}$. For example, in the quaternionic
case, an automorphism of $\mathbb{H}^{n}$ preserves all three complex
structures $I$, $J$ and $K$ on it, in fact it preserves the whole $S^{2}$
(twistor) family of complex structures. From the metric point of view, it is
more natural to allow $\phi$ to rotate these complex structures. This brings
us to the following definition.

\begin{definition}
Suppose $V$ is a normed linear $\mathbb{A}$-space of rank $n$. A $\mathbb{R}%
$-linear isometry $\phi$ of $V$ is called a twisted isomorphism if there
exists $\theta\in SO\left(  \mathbb{A}\right)  $ such that
\[
\phi\left(  vx\right)  =\phi\left(  v\right)  \theta\left(  x\right)
\]
for any $v\in V$ and $x\in\mathbb{A}$.

We denote the group of twisted isomorphisms of $%
{\textstyle\bigoplus^{n}}
\mathbb{A}$ as $G_{\mathbb{A}}\left(  n\right)  $.
\end{definition}

The following proposition identifies the group of twisted isomorphisms for
each $\mathbb{A}$. Recall that when $\mathbb{A}=\mathbb{O}$, we assume that
$V$ has octonion dimension one.

\begin{proposition}
The normed algebras and their corresponding groups of twisted isomorphisms are
given as follows,
\[%
\begin{array}
[c]{cccccc}%
\mathbb{A} & = & \mathbb{R} & \mathbb{C} & \mathbb{H} & \mathbb{O}\\
G_{\mathbb{A}}\left(  n\right)  & = & O\left(  n\right)  & U\left(  n\right)
& Sp\left(  n\right)  Sp\left(  1\right)  & Spin\left(  7\right)  .
\end{array}
\]
\end{proposition}

Proof: For $\mathbb{A}=\mathbb{R}$, the assertion is trivial because
$SO\left(  \mathbb{R}\right)  $ is the trivial group. In the complex case, we
have $\theta\in SO\left(  \mathbb{C}\right)  =U\left(  1\right)  $, that is
there exists $z\in\mathbb{C}$ with $\left|  z\right|  =1$ satisfying
$\theta\left(  x\right)  =zx$ for any $x\in\mathbb{C}$. The requirement
$\phi\left(  vx\right)  =\phi\left(  v\right)  zx$ for $x=1$ implies that
$z=1$. That is $\phi$ is complex linear and hence $G_{\mathbb{A}}\left(
n\right)  =GL\left(  n,\mathbb{C}\right)  \cap O\left(  2n\right)  =U\left(
n\right)  $.

In the quaternionic case, we have $\theta\in SO\left(  \mathbb{H}\right)
=Sp\left(  1\right)  Sp\left(  1\right)  $, i.e. $\theta\left(  x\right)
=\alpha x\beta$ for some unit quaternions $\alpha,\beta\in S^{3}%
\subset\mathbb{H}$. The requirement $\phi\left(  vx\right)  =\phi\left(
v\right)  \alpha x\beta$ with $x=1$ implies that $\alpha\beta=1\in SO\left(
4\right)  $, i.e. $\beta=\alpha^{-1}$. If we define a linear homomorphism $A$
by
\[
A\left(  v\right)  =\phi\left(  v\right)  \alpha,
\]
for any $v\in\mathbb{H}^{n}$, then $A\left(  vx\right)  =A\left(  v\right)  x$
for any $x\in\mathbb{H}$, i.e. $A\in GL\left(  n,\mathbb{H}\right)  \cap
O\left(  4n\right)  =Sp\left(  n\right)  $. From this, we have $\phi\in
Sp\left(  n\right)  Sp\left(  1\right)  =G_{\mathbb{H}}\left(  n\right)  $.
For the octonionic case, the identification of $G_{\mathbb{O}}$ with
$Spin\left(  7\right)  $ can be found in \cite{Murakami}. $\blacksquare$

\subsection{\label{1Sec Riem A mfd}Special Riemannian $\mathbb{A}$-manifolds}

We begin with the definition of Riemannian manifolds defined over $\mathbb{A}
$.

\begin{definition}
A Riemannian manifold $\left(  M,g\right)  $ of is called a Riemannian
$\mathbb{A}$-manifold, or simply an $\mathbb{A}$-manifold, if the holonomy
group of its Levi-Civita connection is a subgroup of $G_{\mathbb{A}}\left(
n\right)  \subset O\left(  m\right)  $ with $m=\dim M=n\dim\mathbb{A}$.
\end{definition}

From the previous proposition, we know that Riemannian $\mathbb{A}$-manifolds
for various $\mathbb{A}$ have holonomy groups insides $O\left(  n\right)  $,
$U\left(  n\right)  $, $Sp\left(  n\right)  Sp\left(  1\right)  $ and
$Spin\left(  7\right)  $ respectively and these manifolds are called
Riemannian manifolds, K\"{a}hler manifolds, quaternionic K\"{a}hler manifolds
and $Spin\left(  7\right)  $-manifolds respectively.

In section \ref{1Sec Rem Qu}, we will discuss $\mathbb{A}$-manifolds without
Riemannian metrics, e.g. complex manifolds.

\bigskip

Next we introduce the notion of an $\mathbb{A}$-orientation for Riemannian
$\mathbb{A}$-manifolds. For a real Riemannian manifold $M$, an orientation is
a parallel volume form on $M$. Equivalently the holonomy group of $M$ is
inside $SO\left(  m\right)  \subset O\left(  m\right)  $. The determinant
defines a natural action of $O\left(  m\right)  $ on $\Lambda^{m}%
\mathbb{R}\cong\mathbb{R}$ and $SO\left(  m\right)  $ is the isotropic
subgroup for any nonzero element in $\mathbb{R}$. In general there is a
natural real representation of $G_{\mathbb{A}}\left(  n\right)  $ on
$\mathbb{A}$,%

\[
\lambda_{\mathbb{A}}:G_{\mathbb{A}}\left(  n\right)  \rightarrow O\left(
\mathbb{A}\right)  \text{,}%
\]
except in the quaternionic case where the action is only defined projectively.

\begin{definition}
For any $g\in G_{\mathbb{A}}\left(  n\right)  \ $and $x\in\mathbb{A}$, we define

(1) $\lambda_{\mathbb{R}}\left(  g\right)  \left(  x\right)  =x\det\left(
g\right)  \in\mathbb{R}$,

(2) $\lambda_{\mathbb{C}}\left(  g\right)  \left(  x\right)  =x\det
{}_{\mathbb{C}}\left(  g\right)  \in\mathbb{C}$,

(3) $\lambda_{\mathbb{H}}\left(  g\right)  \left(  x\right)  =x\beta
\in\mathbb{H}$ with $g=\left(  \alpha,\beta\right)  \in Sp\left(  n\right)
Sp\left(  1\right)  \ $and

(4) $\lambda_{\mathbb{O}}\left(  g\right)  \left(  x\right)  =g\cdot
x\in\mathbb{O}$ via $Spin\left(  7\right)  \subset SO\left(  8\right)  $.
\end{definition}

Note that $G_{\mathbb{A}}\left(  n\right)  $ always acts transitively on the
unit sphere in $\mathbb{A}$.

The explanation of the seemingly different looking $\lambda_{\mathbb{H}}$ is
as follow: First $\det_{\mathbb{H}}$ can not be defined because of the
non-commutativity of $\mathbb{H}$. Even when $\mathbb{A}=\mathbb{C}$,
$\det_{\mathbb{C}}$ can be interpreted as giving a decomposition,
\begin{align*}
&  U\left(  n\right)  \overset{\cong}{\rightarrow}SU\left(  n\right)  \times
U\left(  1\right)  /\mathbb{Z}_{n}\\
A  &  \rightarrow\left(  A\cdot\left(  \det{}_{\mathbb{C}}A\right)
^{-1/n},\,\left(  \det{}_{\mathbb{C}}A\right)  ^{1/n}\right)
\end{align*}
and the projection to the second factor $\left(  \det{}_{\mathbb{C}}A\right)
^{1/n}\in U\left(  1\right)  /\mathbb{Z}_{n}$ can be reinterpreted as an
element in $U\left(  1\right)  $ by raising to the $n^{th}$-power, thus giving
the complex determinant of $A$. The natural analog in the quaternionic case is
the identification,
\[
G_{\mathbb{H}}\left(  n\right)  =Sp\left(  n\right)  Sp\left(  1\right)
\overset{\cong}{\rightarrow}Sp\left(  n\right)  \times Sp\left(  1\right)
/\mathbb{Z}_{2}%
\]
thus $\left(  \alpha,\beta\right)  \rightarrow\beta$ is the direct analog to
$\,\left(  \det{}_{\mathbb{C}}A\right)  ^{1/n}$ in the complex case. Later
this $\mathbb{Z}_{2}$ factor will identify special $\mathbb{C}$-Lagrangian
submanifolds of type I and type II in any hyperk\"{a}hler manifolds to the
same kind of objects, namely the complex Lagrangian submanifolds, see section
\ref{1Sec B Lag submfd}.

\begin{definition}
A twisted isomorphism $g\in G_{\mathbb{A}}\left(  n\right)  $ is called
special if $\lambda_{\mathbb{A}}\left(  g\right)  $ fixes $1\in\mathbb{A}$.
That is $g$ is an element in the isotropic subgroup of $1$ in $G_{\mathbb{A}%
}\left(  n\right)  $, which we denote $H_{\mathbb{A}}\left(  n\right)  $.
\end{definition}

\bigskip

It is not difficult to identify these groups when $\mathbb{A}=\mathbb{R}$,
$\mathbb{C}$ or $\mathbb{H}$. In the octonionic case, it is a classical result
(see e.g. \cite{Murakami}) that $H_{\mathbb{O}}\left(  1\right)  $ is $G_{2}$.
Recall that the Lie group $G_{2}$can be identified as the stabilizer of the
natunal action of $Spin\left(  7\right)  $ on $S^{7}$. Thus we have the
following lemma.

\begin{lemma}
The normed algebras and their corresponding groups of special twisted
isomorphisms are given as follows,
\[%
\begin{array}
[c]{cccccc}%
\mathbb{A} & = & \mathbb{R} & \mathbb{C} & \mathbb{H} & \mathbb{O}\\
H_{\mathbb{A}}\left(  n\right)  & = & SO\left(  n\right)  & SU\left(  n\right)
& Sp\left(  n\right)  & G_{2}%
\end{array}
\]
\end{lemma}

\bigskip

Remark: An isomorphism of a normed linear $\mathbb{A}$-space is the same as
(i) a twisted isomorphism when $\mathbb{A}$ is $\mathbb{R}$ or $\mathbb{C}$
and (ii) a special twisted isomorphism when $\mathbb{A}$ is $\mathbb{H}$ or
$\mathbb{O}$.

Next we define the analog of an orientation for Riemannian $\mathbb{A}%
$-manifold $\left(  M,g\right)  $, i.e. the holonomy group for the Levi-Civita
connection is inside $G_{\mathbb{A}}\left(  n\right)  $.

\begin{definition}
A Riemannian $\mathbb{A}$-manifold is called special if its holonomy group is
inside $H_{\mathbb{A}}\left(  n\right)  $.
\end{definition}

There is another characterization of special $\mathbb{A}$-manifolds. Using the
representation $\lambda_{\mathbb{A}}\left(  n\right)  $ of $G_{\mathbb{A}%
}\left(  n\right)  $, we obtain a vector bundle $\mathbb{A}_{M}$ over any
Riemannian $\mathbb{A}$-manifold $M$,%
\[
\mathbb{A}\rightarrow\mathbb{A}_{M}\rightarrow M\text{.}%
\]
Then $M$ is special if and only if there is a parallel section for
$\mathbb{A}_{M}$.

In the real case, we have $\mathbb{R}_{M}=\Lambda^{n}T_{M}^{\ast}$. In the
complex case, we have $\mathbb{C}_{M}=\Lambda^{n}T_{M}^{\ast(1,0)}$, the
canonical line bundle of $M$. In the quaternionic case, the projection of its
holonomy group $Sp\left(  n\right)  Sp\left(  1\right)  $ to its first and
second factor are only well-defined up to $\pm1$. Suppose that $M$ is spin,
$w_{2}\left(  M\right)  =0$, then this $\pm1$ ambiguity can be lifted and we
obtain a $Sp\left(  n\right)  $-bundle $V$ and a $Sp\left(  1\right)  $-bundle
$S$ over $M$ via the standard representation of $Sp\left(  n\right)  $ and
$Sp\left(  1\right)  $ respectively. The inclusion of $Sp\left(  n\right)
Sp\left(  1\right)  $ in $SO\left(  4n\right)  $ gives the isomorphism%

\[
T_{M}^{\ast}\otimes_{\mathbb{R}}\mathbb{C}=V\otimes_{\mathbb{C}}S\text{.}%
\]
We have $\mathbb{H}_{M}=S$. When $w_{2}\left(  M\right)  \neq0$, then $V$ and
$S$ only exist locally (see e.g. \cite{Besse}).

In the octonionic case, we simply have $\mathbb{O}_{M}=T_{M}^{\ast}$.

The following table lists the possible holonomy groups for all (special)
Riemannian $\mathbb{A}$-manifolds and their usual names.%

\[%
\begin{tabular}
[c]{|l|l|l|}\hline
$%
\begin{array}
[c]{l}%
\,\\
\,
\end{array}
$ & Riemannian $\mathbb{A}$-manifolds & Special Riemannian $\mathbb{A}%
$-manifolds\\\hline
$\mathbb{R}%
\begin{array}
[c]{l}%
\,\\
\,
\end{array}
$ & $O\left(  n\right)  $ & $SO\left(  n\right)  $\\
& (Riemannian manifolds) & (O$\text{riented Riemannian manifolds}$)\\\hline
$\mathbb{C}%
\begin{array}
[c]{l}%
\,\\
\,
\end{array}
$ & $U\left(  n\right)  $ & $SU\left(  n\right)  $\\
& (Kahler$\text{ manifolds}$)\thinspace & (Calabi-Yau manifolds)\\\hline
$\mathbb{H}%
\begin{array}
[c]{l}%
\,\\
\,
\end{array}
$ & $Sp\left(  n\right)  Sp\left(  1\right)  $ & $Sp\left(  n\right)  $\\
& (Q$\text{uaternionic-Kahler manifolds}$) & (H$\text{yperkahler manifolds}%
$)\\\hline
$\mathbb{O}%
\begin{array}
[c]{l}%
\,\\
\,
\end{array}
$ & $Spin\left(  7\right)  $ & $G_{2}$\\
& ($Spin\left(  7\right)  \text{-manifolds}$) & ($G_{2}\text{-manifolds}%
$)\\\hline
\end{tabular}
\]

This list gives all possible holonomy groups of a (non-locally symmetric)
irreducible Riemannian manifold $M$, as classified by Berger \cite{Berger}. To
phrase this differently, Riemannian manifolds with various holonomy groups are
classified in terms of a normed division algebra $\mathbb{A}$ and its
$\mathbb{A}$-orientability.

\subsection{\label{1Sec Char Prop}Characterizations and properties}

For completeness, we briefly describe these (special) $\mathbb{A}$-manifolds
and introduce certain natural differential forms on them that we will need
later, see \cite{Besse} or \cite{Salamon Book} for details. In the real case,
(special) Riemannian $\mathbb{R}$-manifolds are simply (oriented) Riemannian
manifolds. Orientability of $M$ allows us to determine a square root
$\sqrt{\det\left(  g_{ij}\right)  }$ consistently, and we obtain a parallel
volume form $\nu_{M}$. Other (special) $\mathbb{A}$-manifolds also admit
characterizations in terms of the existence of certain non-degenerate parallel
forms. In all cases, harmonicity is already enough.

(1) $\mathbb{C}$\textit{-manifolds} (i.e. \textit{K\"{a}hler manifolds}),
$Hol\left(  g\right)  \subset U\left(  n\right)  $. Since $U\left(  n\right)
=O\left(  2n\right)  \cap GL\left(  n,\mathbb{C}\right)  $, the Levi-Civita
connection $\nabla$ preserves an (almost) complex structure $J$, i.e.
\[
J^{2}=-1\text{ and }\nabla J=0\text{.}%
\]
This implies integrability of $J$. Alternatively we can use the K\"{a}hler
form $\omega$,
\[
\omega\left(  v,w\right)  =g\left(  Jv,w\right)  \text{ and }\nabla\omega=0,
\]
to characterize a K\"{a}hler manifold. Namely, a Hermitian metric on a complex
manifold $M$ is K\"{a}hler if and only if $d\omega=0$. It follows that every
projective algebraic manifold inside $\mathbb{CP}^{N}$ is K\"{a}hler.

(2) \textit{Special} $\mathbb{C}$\textit{-manifolds} (i.e. \textit{Calabi-Yau
manifolds}), $Hol\left(  g\right)  \subset SU\left(  n\right)  $. By
definition, these are K\"{a}hler manifolds with a parallel section of the
canonical line bundle $K_{M}=\Lambda^{n}T_{M}^{^{\ast}\left(  1,0\right)  }$,
i.e. a parallel holomorphic volume form $\Omega\in\Omega^{n,0}\left(
M\right)  $. By the celebrated result of Yau \cite{Yau Calabi}, such a
structure always exist on any compact K\"{a}hler manifold with topological
trivial $K_{M}$. This follows that a degree $d$ smooth hypersurface in
$\mathbb{CP}^{n+1}$ admits a Calabi-Yau metric if and only if $d=n+2$. For
instance, when $n=2$, we have a quartic surface in $\mathbb{CP}^{3}$, that is
a K3 surface.

(3) \textit{Special} $\mathbb{H}$\textit{-manifolds} (i.e. \textit{hyperk\"{a}%
hler manifolds}), $Hol\left(  g\right)  \subset Sp\left(  n\right)  $. Since
$Sp\left(  n\right)  $ is the automorphism group of the quaternionic vector
space $\mathbb{H}^{n}$, the Riemannian metric on such a manifold support three
K\"{a}hlerian complex structures $I$, $J$, $K$ satisfying the Hamilton
relation
\[
I^{2}=J^{2}=K^{2}=IJK=-1\text{.}%
\]
A characterization by Hitchin says that a Riemannian metric on $M$ which is
Hermitian with respect to three almost complex structures $I$, $J$ and $K$,
satisfying the Hamilton relation and $d\omega_{I}=d\omega_{J}=d\omega_{K}=0$,
then its holonomy group is inside $Sp\left(  n\right)  $. This implies that
the distinction between K\"{a}hler manifolds and symplectic manifolds no
longer exist in the quaternionic case.

If we denote $\Omega_{J}=\omega_{I}+i\omega_{K}$, then $\Omega_{J}$ is a
parallel holomorphic symplectic form on $M$. Using Yau's theorem, every
compact K\"{a}hler manifold with a holomorphic symplectic form admits a
hyperk\"{a}hler structure.

When $\dim M=4$, i.e. $n=1$, hyperk\"{a}hler manifolds are the same as
Calabi-Yau manifolds because of $Sp\left(  1\right)  =SU\left(  2\right)  $.
If $M$ is compact then it is either a flat torus of dimension four or a K3
surface. Using Yau's theorem, Fujiki and Beauville show that if $S$ is a
compact hyperk\"{a}hler manifold of dimension four, then the minimal
resolution of the symmetric product of $S$ admits a natural hyperk\"{a}hler structure.

(4) $\mathbb{H}$\textit{-manifold}, (i.e. \textit{quaternionic K\"{a}hler
manifolds}), $Hol\left(  g\right)  \subset Sp\left(  n\right)  Sp\left(
1\right)  $. Such a manifold is similar to a hyperk\"{a}hler manifold,
however, the individual complex structures $I$, $J$ and $K$ can only be
defined locally. The four form $\Theta=\omega_{I}^{2}+\omega_{J}^{2}%
+\omega_{K}^{2}$ is nevertheless well-defined and parallel. Examples include
quaternionic projective spaces $\mathbb{HP}^{n}$.

(5) \textit{Special} $\mathbb{O}$\textit{-manifolds}, (i.e. $G_{2}%
$\textit{-manifolds}), $Hol\left(  g\right)  \subset G_{2}$. Since
$G_{2}=H_{\mathbb{O}}$ stabilizes $1\in\mathbb{O}$, it is really a subgroup of
$SO\left(  \operatorname{Im}\mathbb{O}\right)  =SO\left(  7\right)  $. That
implies that, up to covering, $M=X\times\mathbb{R}$ because of the deRham
decomposition (see e.g. \cite{Besse}). Traditionally a $G_{2}$-manifold is
referred to the \textit{seven} dimensional manifold $X$. The cross product on
$\mathbb{O}$, defined as $x\times y=\operatorname{Im}\bar{y}x$, induces a
product structure $\times$ on any $G_{2}$-manifold because $G_{2}$ is the
automorphism group of the normed algebra $\mathbb{O}$. This determines a
parallel (positive) three form $\Omega$,
\[
\Omega\left(  x,y,z\right)  =\left\langle x,y\times z\right\rangle \text{ and
}\triangledown\Omega=0
\]
on any $G_{2}$-manifold $X$. Recall that the natural action of $GL\left(
7,\mathbb{R}\right)  $ on $\Lambda^{3}\mathbb{R}^{7}$ has two open orbits,
called positive and negative (see e.g. \cite{Hitchin 3 form}). Those three
forms in the same open orbit as the one given in the following table are
called positive. Gray shows that any harmonic positive three form on $X$
determines a $G_{2}$-manifold structure on it. We will also use the parallel
four form $\Theta=\ast\Omega$ later.

(6) $\mathbb{O}$\textit{-manifolds}, (i.e. $Spin\left(  7\right)
$\textit{-manifolds}), $Hol\left(  g\right)  \subset Spin\left(  7\right)
\subset SO\left(  8\right)  $. Gray shows that an eight dimensional manifold
$M$ has holonomy group $Spin\left(  7\right)  $ if and only if it admits a
harmonic (and hence parallel) self-dual four form $\Theta$. We have
$\Theta=\Omega_{G_{2}}\wedge dx^{0}-\Theta_{G_{2}}$ when our $\mathbb{O}%
$-manifold $M=X\times\mathbb{R}$ is special. Complete examples of (special)
$\mathbb{O}$-manifolds are constructed by Bryant and Salamon
\cite{BryantSalamon} and compact examples are constructed by Joyce \cite{Joyce
book}, and recently by Kovalev \cite{Kovalev}.

In the next table, we list the various parallel forms on $\mathbb{A}%
$-manifolds, in the Euclidean case:\footnote{$\left[  123\right]  $ means
terms obtained by permuting the indexes insides the bracket.}%

\[%
\begin{tabular}
[c]{|l|l|ll|}\hline
$%
\begin{array}
[c]{l}%
\,\\
\,
\end{array}
$ & Riemannian $\mathbb{A}$-manifolds &  & Parallel forms\\\hline\hline
$\mathbb{R}$ & Oriented Riem. manifolds & On $\mathbb{R}^{n}$, & $\nu
_{M}=dx^{1}dx^{2}\cdots dx^{n}%
\begin{array}
[c]{l}%
\,\\
\,
\end{array}
$\\\hline
$\mathbb{C}$ & K\"{a}hler manifolds & On $\mathbb{C}^{n}$, & $\omega=\frac
{i}{2}\left(  dz^{1}d\bar{z}^{1}+\cdots+dz^{n}d\bar{z}^{n}\right)
\begin{array}
[c]{l}%
\,\\
\,
\end{array}
$\\
& Calabi-Yau manifolds & On $\mathbb{C}^{n}$, & $\omega$ and $\Omega
=dz^{1}\wedge\cdots\wedge dz^{n}%
\begin{array}
[c]{l}%
\,\\
\,
\end{array}
$\\\hline
$\mathbb{H}$ & Quaternionic K\"{a}hler mfds & On $\mathbb{H}^{n}$, &
$\Theta=\omega_{I}^{2}+\omega_{J}^{2}+\omega_{K}^{2}%
\begin{array}
[c]{l}%
\,\\
\,
\end{array}
$\\
& Hyperk\"{a}hler manifolds & On $\mathbb{H}^{n}$, & $\omega_{J}=\frac{i}%
{2}\left(  dz^{1}d\bar{z}^{1}+\cdots+dz^{2n}d\bar{z}^{2n}\right)  $\\
&  &  & $\Omega_{J}=\omega_{I}+i\omega_{K}=dz^{1}dz^{2}+\cdots+dz^{2n-1}%
dz^{2n}$\\\hline
$\mathbb{O}$ & $Spin\left(  7\right)  $-manifolds & On $\mathbb{R}^{8}$, &
$\Theta=\Omega_{G_{2}}\wedge dx^{0}-\Theta_{G_{2}}%
\begin{array}
[c]{l}%
\,\\
\,
\end{array}
$\\
& $G_{2}$-manifolds & On $\mathbb{R}^{7}$, & $\Omega_{G_{2}}=dx^{123}%
-dx^{1}\left(  dy^{23}+dy^{10}\right)  +\left[  123\right]  $\\
&  &  & $\Theta_{G_{2}}=\ast\Omega_{G_{2}}$\\\hline
\end{tabular}
\]

\bigskip

Remark: Special Riemannian $\mathbb{O}$-manifolds (with $\mathbb{O}%
\neq\mathbb{R}$) are Calabi-Yau manifolds, hyperk\"{a}hler manifolds and
$G_{2}$-manifolds. These spaces are the central objects of interests in string
theory, conformal field theory and M-theory in physics. From a mathematical
point of view, they share many good geometric properties\footnote{Many of
these properties are also shared by $Spin\left(  7\right)  $-manifolds.}: (1)
their Ricci tensors are all zero, $Ricci=0$. (2) If $M$ is compact with
holonomy group equals $H_{\mathbb{A}}\left(  n\right)  $, then $\pi_{1}\left(
M\right)  $ is finite. (3) the moduli spaces of these metrics are always
smooth, as shown by Bogomolov, Joyce, Tian and Todorov. (4) We can define a
\textit{period map} on the moduli space by integrating the parallel form
$\Omega$ over topological cycles. Locally the period map determines the moduli
space together with its Weil-Peterrson metric. (5) The first Pontriajin number
of $M$ is non-negative, moreover it is zero if and only if $M$ is flat.

\bigskip

On real (resp. complex) manifolds, there are many local differentiable (resp.
holomorphic) functions, which are used to describe the geometry of these
manifolds. However on $\mathbb{H\,}$- (or $\mathbb{O\,}$-)manifolds, there are
very few such functions. For instance, Gray shows that every quaternionic map
is totally geodesic. In particular there are no quaternionic submanifolds in
$\mathbb{HP}^{N}$ other than affine subspaces. Instead the geometries of these
manifolds are reflected by their (Yang-Mills calibrated) $\mathbb{A}%
$-connections and (volume calibrated) $\mathbb{A}$-submanifolds and $\frac
{1}{2}\mathbb{A}$-Lagrangian submanifolds.

\section{\label{1Sec YM bdl}Yang-Mills bundles}

Suppose $E$ is a Hermitian vector bundle over $M$%
\[
\mathbb{C}^{r}\rightarrow E\rightarrow M\text{.}%
\]
A connection on $E$ gives a first order differential operator
\[
D_{E}:\Omega^{k}\left(  M,E\right)  \rightarrow\Omega^{k+1}\left(  M,E\right)
\text{.}%
\]
Its square is a zeroth order operator, called the curvature tensor
$F_{E}=\left(  D_{E}\right)  ^{2}\in\Lambda^{2}\left(  M,ad\left(  E\right)
\right)  $. If $D_{E}$ is a flat connection, i.e. $F_{E}=0$, then holonomy
around any point $x\in M$ gives a representation of the fundamental group,
\[
\rho:\pi_{1}\left(  M,x\right)  \rightarrow U\left(  r\right)  \text{,}%
\]
and vice versa.

When $M$ is a complex manifold of complex dimension $n$, we can decompose
differential forms into $\left(  p,q\right)  $-types,
\[
\Lambda^{2}\left(  M,\mathbb{C}\right)  =\Lambda^{2,0}\left(  M\right)
+\Lambda^{1,1}\left(  M\right)  +\Lambda^{0,2}\left(  M\right)  \text{.}%
\]
If $E$ is a \textit{holomorphic }bundle over $M$ then it admits a connection
with $F_{E}^{0,2}=F_{E}^{2,0}=0$ and vice versa. Roughly speaking, holomorphic
structures on $E$ are equivalent to \textit{partial flat} connections on it.
When $M$ is also K\"{a}hler, we have a further decomposition into primitive
components (e.g. \cite{Besse}),
\[
\Lambda^{1,1}\left(  M\right)  =\Lambda_{0}^{1,1}\left(  M\right)
+\mathbb{C}\omega\text{.}%
\]
Hermitian connections on $E$ satisfying $F_{E}^{0,2}=F_{E}^{2,0}=0$ and
$F_{E}\wedge\omega^{n-1}=0$ are called \textit{Hermitian Yang-Mills}
connections (with zero slope). They have absolute minimum Yang-Mills energy
(see section \ref{1Sec YM calib}). By a famous result of Donaldson, Uhlenbeck
and Yau, the existence of such connections is equivalent to the Mumford
poly-stability for $E$, a natural notion in algebraic geometry which is used
in constructing moduli spaces. Note that we can rephrase these two equations
as requiring $F_{E}\in\Lambda_{0}^{1,1}\left(  M,ad\left(  E\right)  \right)
$.

Another familiar class of partial flat connections are ASD connections over an
oriented Riemannian four manifold $M$. Using the isomorphism $SO\left(
4\right)  \cong Sp\left(  1\right)  Sp\left(  1\right)  $ we can decompose two
forms into self-dual and anti-self-dual components,
\begin{align*}
\Lambda^{2}\left(  M\right)   &  =\Lambda_{+}^{2}\left(  M\right)
+\Lambda_{-}^{2}\left(  M\right) \\
F  &  =F^{+}+F^{-}=\frac{F+\ast F}{2}+\frac{F-\ast F}{2}\text{.}%
\end{align*}

A Hermitian connection $D_{E}$ over $M$ is called an \textit{instanton}, or
\textit{ASD} connection if the self-dual component of its curvature vanishes,
$F_{E}^{+}=0$. Equivalently we have $F_{E}\in\Lambda_{-}^{2}\left(
M,ad\left(  E\right)  \right)  $. Again an ASD connection always have absolute
minimum Yang-Mills energy. Donaldson studies their moduli space in details and
obtain many nontrivial results in four dimensional differential topology
\cite{Don Instanton Book}.

We are going to generalize these and define natural classes of partial flat
connections over Riemannian $\mathbb{A}$-manifolds.

\subsection{\label{1Sec Sp A conn}(Special) $\mathbb{A}$-connections}

Suppose that $M$ is a (special) Riemannian $\mathbb{A}$-manifold of dimension
$m=2^{a}n$, where $a=\dim\mathbb{A}$. Its holonomy group $\mathbf{G}%
_{\mathbb{A}}\left(  n\right)  $ (resp. $\mathbf{H}_{\mathbb{A}}\left(
n\right)  $) is a subgroup of $O\left(  m\right)  $. We denote their Lie
algebra as $\mathbf{g}_{\mathbb{A}}\left(  n\right)  $ and $\mathbf{h}%
_{\mathbb{A}}\left(  n\right)  $ respectively. We have%

\[
\mathbf{h}_{\mathbb{A}}\left(  n\right)  \subset\mathbf{g}_{\mathbb{A}}\left(
n\right)  \subset\mathbf{o}\left(  m\right)  .
\]
Using the natural identification $\mathbf{o}\left(  m\right)  \cong\Lambda
^{2}\mathbb{R}^{m\ast}$, we obtain natural subbundles
\[
\mathbf{h}_{\mathbb{A}}\left(  T_{M}\right)  \subset\mathbf{g}_{\mathbb{A}%
}\left(  T_{M}\right)  \subset\Lambda^{2}\left(  M\right)  ,
\]
over a (special) $\mathbb{A}$-manifold $M$. In fact the subbundles
$\mathbf{h}_{\mathbb{A}}\left(  T_{M}\right)  \subset\Lambda^{2}\left(
M\right)  $ is well-defined even for $\mathbb{A}$-manifolds as long as
$\mathbb{A}\neq\mathbb{O}$. This is because $\mathbf{h}_{\mathbb{A}}\left(
n\right)  $ is an ideal in $\mathbf{g}_{\mathbb{A}}\left(  n\right)  $ in
these cases. We define two natural classes of partial flat connections on an
$\mathbb{A}$-manifold as follows.

\begin{definition}
Suppose $D_{E}$ is a connection on a vector bundle $E$ over a (special)
Riemannian $\mathbb{A}$-manifold $M$. We denote its curvature two form as
$F_{E}\in\Lambda^{2}\left(  M,ad\left(  E\right)  \right)  $. Then

(1) $D_{A}$ a called an $\mathbb{A}$-connection if
\[
F_{A}\in\mathbf{g}_{\mathbb{A}}\left(  T_{M}\right)  \otimes ad\left(
E\right)  ,
\]

(2) $D_{A}$ is called a special $\mathbb{A}$-connection if
\[
F_{A}\in\mathbf{h}_{\mathbb{A}}\left(  T_{M}\right)  \otimes ad\left(
E\right)  .
\]
\end{definition}

We are going to describe these partial flat connections individually. Most of
them can be identified with well-known Yang-Mills connections in the
literature, as indicated in the following table.
\[%
\begin{tabular}
[c]{|l|l|l|}\hline
$%
\begin{array}
[c]{l}%
\,\\
\,
\end{array}
$ & $\mathbb{A}$-connections & special $\mathbb{A}$-connections\\\hline
$\mathbb{C}$ & $F_{E}^{0,2}=F_{E}^{2,0}=0%
\begin{array}
[c]{l}%
\,\\
\,
\end{array}
$ & $F_{E}^{0,2}=F_{E}^{2,0}=\Lambda F=0$\\
& (Holomorphic bundles)\thinspace & (Hermitian Yang-Mills bdls.)\\\hline
$\mathbb{H}$ & $F\in\mathbf{g}_{\mathbb{H}}\left(  T_{M}\right)  \otimes
ad\left(  E\right)  $ & $F_{I}^{0,2}=F_{J}^{0,2}=F_{K}^{0,2}=0%
\begin{array}
[c]{l}%
\,\\
\,
\end{array}
$\\
&  & (ASD or hyperholomorphic bdls.)\\\hline
$\mathbb{O}$ & $\ast F_{E}+\Theta\wedge F_{E}=0$ & $F\wedge\Theta=0%
\begin{array}
[c]{l}%
\,\\
\,
\end{array}
$\\
& ($Spin\left(  7\right)  $-Donaldson-Thomas bdls.) & ($G_{2}$%
-Donaldson-Thomas bdls.)\\\hline
\end{tabular}
\]

When $\mathbb{A}=\mathbb{R}$ we have $\mathbf{o}\left(  n\right)
=\mathbf{h}_{\mathbb{R}}\left(  n\right)  =\mathbf{g}_{\mathbb{R}}\left(
n\right)  $. Therefore a (special) $\mathbb{R}$-connection is simply any
connection over $M$.

When $\mathbb{A}=\mathbb{C}$ we have $\mathbf{g}_{\mathbb{C}}\left(  n\right)
=\mathbf{u}\left(  n\right)  $ and $\mathbf{h}_{\mathbb{C}}\left(  n\right)
=\mathbf{su}\left(  n\right)  $. The Lie algebra $\mathbf{u}\left(  n\right)
$ consists of skew-Hermitian matrices. Using the Hermitian inner product on
the vector space $V\cong\mathbb{C}^{n}$ to identify $V$ with $\bar{V}^{\ast}$,
we obtain an identification $\mathbf{u}\left(  n\right)  =\left(  V^{\ast
}\otimes\bar{V}^{\ast}\right)  \cap\Lambda^{2}V_{\mathbb{R}}^{\ast}$. Globally
on $M$, we get $\mathbf{g}_{\mathbb{c}}\left(  T_{M}\right)  \cong
\Lambda^{1,1}\left(  M\right)  _{\mathbb{R}}$. Similarly the trace component
in $\mathbf{u}\left(  n\right)  $ corresponds to the $\mathbb{R}$-span of the
K\"{a}hler form $\omega$ in $\Lambda^{1,1}\left(  M\right)  _{\mathbb{R}}$.
Thus in the primitive decomposition of two forms on $M$,
\[
\Lambda^{2}\left(  M\right)  =\Lambda_{0}^{1,1}\left(  M\right)  _{\mathbb{R}%
}+\mathbb{R}\omega+\left(  \Lambda^{2,0}\left(  M\right)  +\Lambda
^{0,2}\left(  M\right)  \right)  _{\mathbb{R}}\text{,}%
\]
we have
\begin{align*}
\mathbf{g}_{\mathbb{c}}\left(  T_{M}\right)   &  \cong\Lambda^{1,1}\left(
M\right)  _{\mathbb{R}}\\
\mathbf{h}_{\mathbb{c}}\left(  T_{M}\right)   &  \cong\Lambda_{0}^{1,1}\left(
M\right)  _{\mathbb{R}}\text{.}%
\end{align*}
Therefore a Hermitian connection $D_{E}$ over $M$ is a $\mathbb{C}$-connection
iff it defines a holomorphic structure on $E$ and it is special if
$F_{E}\wedge\omega^{n-1}=0$, i.e. a Hermitian Yang-Mills connection.

When $\mathbb{A}=\mathbb{H}$ we have $\mathbf{g}_{\mathbb{H}}\left(  n\right)
=\mathbf{sp}\left(  n\right)  \mathbf{sp}\left(  1\right)  $, $\mathbf{h}%
_{\mathbb{H}}\left(  n\right)  =\mathbf{sp}\left(  n\right)  $ and%

\[
T_{M}^{\ast}\otimes_{\mathbb{R}}\mathbb{C}=V\otimes_{\mathbb{C}}S\text{,}%
\]
well-defined up to $\otimes L^{1/2}$, for any $\mathbb{H}$-manifold $M$. In
any event, $V\otimes V$ and $S\otimes S$ are always well-defined and we have
the following decomposition of two forms on $M$,
\[
\Lambda^{2}\left(  M,\mathbb{C}\right)  =Sym^{2}V\otimes\Lambda^{2}%
S+\Lambda^{2}V\otimes Sym^{2}S.
\]
Note that (i) $S$ being a $Sp\left(  1\right)  $-bundle implies $\Lambda^{2}S$
is trivial and (ii) $V$ is a symplectic bundle, using its symplectic form
$\Omega\in\Lambda^{2}V^{\ast}\cong\Lambda^{2}V$, we have a further splitting
$\Lambda^{2}V=\Lambda_{0}^{2}V+\mathbb{C}.$ Hence%

\[
\Lambda^{2}\left(  M,\mathbb{C}\right)  =Sym^{2}V+Sym^{2}S+\Lambda_{0}%
^{2}V\otimes Sym^{2}S.
\]

When $\dim_{\mathbb{R}}M=4$, $Sym^{2}S$ (resp. $Sym^{2}V$) coincides with the
space of self-dual (resp. anti-self-dual) two forms on $M$ and $\Lambda
_{0}^{2}V\otimes Sym^{2}S$ is zero. We will continue to call a connection
$D_{E}$ with its curvature $F_{E}$ inside $Sym^{2}V\otimes ad\left(  E\right)
$ an \textit{anti-self-dual connection}, even though it is traditionally
called a B-connection.

If we denote the standard representation of $\mathbf{sp}\left(  n\right)  $ as
$V_{0}$, then the Lie algebra $\mathbf{sp}\left(  n\right)  $ is naturally
identified with $Sym^{2}V_{0}$. This implies that
\begin{align*}
\mathbf{g}_{\mathbb{H}}\left(  T_{M}\right)   &  \cong Sym^{2}V+Sym^{2}S,\\
\mathbf{h}_{\mathbb{H}}\left(  T_{M}\right)   &  \cong Sym^{2}V\text{.}%
\end{align*}
In particular a special $\mathbb{H}$-connection is the same as an
anti-self-dual connection on $M$. Similar to Hermitian Yang-Mills connections,
special $\mathbb{H}$-connections have absolute minimum Yang-Mills energy. On
the other hand, on an oriented four manifold, every connection is a
$\mathbb{H}$-connection because $\mathbf{\ sp}\left(  1\right)  +\mathbf{sp}%
\left(  1\right)  =\mathbf{so}\left(  4\right)  $.

Remark: There is also an identification,
\[
Sym^{2}V=\ \Lambda_{I}^{1,1}\cap\ \Lambda_{J}^{1,1}\cap\Lambda_{K}^{1,1}.
\]
Therefore a connection is a special $\mathbb{H}$-connection if and only if it
is holomorphic with respect to $I$, $J$ and $K$, and it is sometimes called a
\textit{hyperholomorphic} connection. On a hyperk\"{a}hler manifold $M$,
Verbitsky \cite{Verbitsky} shows that if $D_{E}$ is a Hermitian-Yang-Mills
connection with respect to the K\"{a}hler structure $\omega_{J}$ and
$c_{1}\left(  E\right)  $ and $c_{2}\left(  E\right)  $ are both $Sp\left(
1\right)  $ invariant cohomology class, then $D_{E}$ is a special $\mathbb{H}$-connection.

\bigskip

When $\mathbb{A}=\mathbb{O}$ we have $\mathbf{g}_{\mathbb{O}}\left(  n\right)
=\mathbf{spin}\left(  7\right)  $ and $\mathbf{h}_{\mathbb{O}}\left(
n\right)  =\mathbf{g}_{2}$. When $M$ is a $\mathbb{O}$-manifold, i.e. a
$Spin\left(  7\right)  $-manifold, we have a natural decomposition of two
forms \cite{Salamon Book},
\[
\Lambda^{2}\left(  M\right)  =\Lambda_{21}^{2}\left(  M\right)  +\Lambda
_{7}^{2}\left(  M\right)  \text{.}%
\]
They are characterized as follows: for any $\phi\in\Lambda^{2}\left(
M\right)  $,
\begin{align*}
\phi &  \in\Lambda_{21}^{2}\left(  M\right)  \text{ iff }\phi+\ast\left(
\Theta\wedge\phi\right)  =0\\
\phi &  \in\Lambda_{7}^{2}\left(  M\right)  \text{ iff }3\phi=\ast\left(
\Theta\wedge\phi\right)  \text{.}%
\end{align*}
Furthermore, we have
\[
\mathbf{g}_{O}\left(  T_{M}\right)  \cong\Lambda_{21}^{2}\left(  M\right)
\text{.}%
\]
Therefore the curvature of any $\mathbb{O}$-connection satisfies%
\[
F_{E}+\ast\left(  \Theta\wedge F_{E}\right)  =0\text{.}%
\]

This equation, and its $G_{2}$-analog, are introduced by Donaldson and Thomas
in \cite{DT}.

When $M=X\times S^{1}$ is a $G_{2}$-manifold, then corresponding to the
reduction from $SO\left(  8\right)  $ to $SO\left(  7\right)  $ we have a
decomposition,
\[
\Lambda^{2}\left(  M\right)  =\Lambda^{2}\left(  X\right)  +\Lambda^{1}\left(
X\right)  \wedge d\theta
\]
where $\theta$ is the angle coordinate on $S^{1}$. On the other hand we have a
similar decomposition of two forms \cite{Salamon Book} for the seven
dimensional manifold $X$,
\[
\Lambda^{2}\left(  X\right)  =\Lambda_{14}^{2}\left(  X\right)  +\Lambda
_{7}^{2}\left(  X\right)  \text{,}%
\]
with
\[
\mathbf{h}_{\mathbb{O}}\left(  T_{M}\right)  \cong\Lambda_{14}^{2}\left(
X\right)  \text{.}%
\]
A special $\mathbb{O}$-connection on $X$ is again a $G_{2}$-Donaldson-Thomas
connection, i.e.
\[
F_{E}\wedge\Theta=0\text{,}%
\]
and these are absolute minimums for the Yang-Mills energy.

\subsection{\label{1Sec YM calib}Yang-Mills calibrations}

Suppose $E$ is a Hermitian vector bundle over an oriented Riemannian manifold
$M$, with volume form $\nu_{M}$. The \textit{Yang-Mills energy} of a Hermitian
connection $D_{A}$ is defined as follows,%

\[
YM\left(  D_{A}\right)  =\int_{M}\left|  F_{A}\right|  ^{2}v_{M}.
\]
The Euler-Lagrange equation is called the \textit{Yang-Mills equation} and it
is given by,
\[
D_{A}^{\ast}F_{A}=0\text{.}%
\]

Analogous to the volume calibration for minimal submanifolds (see \cite{HL}
and section \ref{1Sec Vol calib}), we have the notion of Yang-Mills
calibration for connections, which gives Yang-Mills connections with absolute
minimal energy. This is a modification of the $\Omega$-ASD connections
introduced by Tian in \cite{Tian YM Cal}.

On a vector space $V\cong\mathbb{R}^{m}$ with a fixed volume form $\nu$, each
element $\Phi$ in $\Lambda^{m-4}V^{\ast}$ defines a quadratic form $q_{\Phi}$
on $\Lambda^{2}V^{\ast}$ as follows,
\begin{align*}
q_{\Phi}  &  :\Lambda^{2}V^{\ast}\rightarrow\mathbb{R}\\
q_{\Phi}\left(  \phi\right)   &  =\phi\wedge\phi\wedge\Phi/\nu_{M}\text{.}%
\end{align*}

\begin{definition}
Suppose $M$ is an oriented Riemannian manifold $M$ of dimension $m$, a
differential form $\Phi\in\Omega^{m-4}\left(  M\right)  $ is called a
Yang-Mills calibrating form if
\begin{gather*}
d\Phi=0\\
q_{\Phi}\left(  \phi\right)  \leq\left|  \phi\right|  ^{2}%
\end{gather*}
for any $\phi\in\Omega^{2}\left(  M\right)  $.
\end{definition}

\begin{definition}
Suppose $E$ is a Hermitian vector bundle over a manifold $M$ with a Yang-Mills
calibrating form $\Phi$. A Hermitian connection $D_{A}$ on $E$ is called
Yang-Mills calibrated by $\Phi$, or simply $\Phi$-calibrated, if its curvature
tensor $F_{A}$ satisfies\footnote{The quadratic form $q_{\Phi}$ is extended to
$ad\left(  E\right)  $-valued two forms using the Killing form.}%
\[
q_{\Phi}\left(  F_{A}\right)  \leq\left|  F_{A}\right|  ^{2}\text{.}%
\]
\end{definition}

\bigskip

As in the volume calibration case, we have the following fundamental lemma.

\begin{lemma}
If $D_{A}$ is a $\Phi$-calibrated connection on $E$ and $D_{A^{\prime}}$ is
any other connection, then the Yang-Mills energy of $D_{A^{\prime}}$ is
smaller than or equal to that for $D_{A}$,
\[
YM\left(  D_{A^{\prime}}\right)  \geq YM\left(  D_{A}\right)  .
\]

Moreover, if the equality sign holds, then $D_{A^{\prime}}$ is also $\Phi$-calibrated.
\end{lemma}

Proof:
\begin{align*}
YM\left(  D_{A^{\prime}}\right)   &  =\int_{M}\left|  F_{A^{\prime}}\right|
^{2}\nu_{M}\geq\int_{M}Tr\left(  F_{A^{\prime}}^{2}\right)  \wedge\Phi\\
&  =\int_{M}Tr\left(  F_{A}^{2}\right)  \wedge\Phi\,\,\,\text{ \thinspace
\ (since }d\Phi=0\text{)}\\
&  =\int_{M}\left|  F_{A}\right|  ^{2}\nu_{M}=YM\left(  D_{A}\right)  \text{.}%
\end{align*}

Hence the result. $\blacksquare$

\bigskip

The following Chern number inequality gives a topological constraint to the
existence of $\Phi$-calibrated connections on $E$. It also give an effective
way to characterize flat connections. The proof of it is simple and standard.

\begin{proposition}
If $E$ admits a $\Phi$-calibrated connection then we have
\[
\int_{M}ch\left(  E\right)  \Phi\leq0
\]
and the equality sign holds iff $E$ is a flat bundle.
\end{proposition}

Proof: This follows immediately from the Chern-Weil formula
\[
ch\left(  E\right)  =\exp\left(  \frac{i}{2\pi}F_{E}\right)
\]
and the definition of a $\Phi$-calibrated connection. $\blacksquare$

\bigskip

For any Riemannian $\mathbb{A}$-manifold $M$, there is a natural Yang-Mills
calibrating form and connections they calibrate are basically the same as
special $\mathbb{A}$-connections over $M$. We list these Yang-Mills
calibrating forms are the common names for the corresponding calibrated
Yang-Mills connections.
\[%
\begin{tabular}
[c]{|l|l|l|}\hline
$%
\begin{array}
[c]{l}%
\,\\
\,
\end{array}
$ & Calibrating form, $\Phi$ & Yang-Mills connections\\\hline\hline
$\mathbb{R}$-manifold & $\Phi=0%
\begin{array}
[c]{l}%
\,\\
\,
\end{array}
$ & Flat connections\\\hline
$\mathbb{C}$-manifold & $\Phi=\omega^{n-2}%
\begin{array}
[c]{l}%
\,\\
\,
\end{array}
$ & Hermitian Yang-Mills connections\\\hline
$\mathbb{H}$-manifold & $\Phi=\Theta^{n-1}%
\begin{array}
[c]{l}%
\,\\
\,
\end{array}
$ & Anti-Self-Dual connections\\\hline
$\mathbb{O}$-manifold & $\Phi=\Theta%
\begin{array}
[c]{l}%
\,\\
\,
\end{array}
$ & Donaldson-Thomas connections\\\hline
\end{tabular}
\]

For example, when $M$ is an oriented Riemannian four manifold, $\Phi=1$ is a
Yang-Mills calibrating form and connections it calibrates are precisely
ASD\ connections, i.e. $F_{A}^{+}=0$.

\section{\label{1Sec Min submfd}Minimal submanifolds}

On any general Riemannian $\mathbb{A}$-manifolds $M$, we introduce two natural
classes of submanifolds: (i) $\mathbb{A}$-submanifolds and (ii) $\frac{1}%
{2}\mathbb{A}$-Lagrangian submanifolds. $\mathbb{A}$-submanifolds can be
defined on any $\mathbb{A}$-manifold, even without a Riemannian metric (see
section \ref{1Sec Rem Qu}). $\frac{1}{2}\mathbb{A}$-Lagrangian submanifolds
can be viewed as the \textit{maximally real submanifolds} in $M$. For example,
when $\mathbb{A}=\mathbb{C}$ they are (i) complex submanifolds and (ii)
Lagrangian submanifolds in a K\"{a}hler manifold.

When $M$ is special, there are two natural subclasses of (ii), called special
$\frac{1}{2}\mathbb{A}$-Lagrangian submanifolds of type I and type II. For
example, if $M$ is a Calabi-Yau manifold, then they are special Lagrangian
submanifolds of phase $0$ and phase $\pi/2$.

These $\mathbb{A}$-submanifolds and special $\frac{1}{2}\mathbb{A}$-Lagrangian
submanifolds are always absolute volume minimizers, as they are volume calibrated.

\subsection{\label{1Sec A submfd}$\mathbb{A}$-submanifolds}

We first discuss the linear case. Suppose $V$ is a linear $\mathbb{A}$-space.
To discuss its linear $\mathbb{A}$-subspaces, we can assume that
$\mathbb{A}\neq\mathbb{O}$, since every nontrivial linear $\mathbb{O}$-space
is isomorphic to $\mathbb{O}$ itself. A linear $\mathbb{A}$-space is the same
as a bi-module over $\mathbb{A}$, which is really an $\mathbb{A}$-vector space
when $\mathbb{A}$ equals $\mathbb{R}$ or $\mathbb{C}$. A linear $\mathbb{A}%
$-subspace of $V$ is then a bi-submodule of $V$, i.e. a real vector subspace
in $V$ which is stable under the left and right action of $\mathbb{A}$. To
globalize these $\mathbb{A}$-subspace structures to any $\mathbb{A}$-manifold,
we need to know that they are stable under twisted isomorphisms.

\begin{lemma}
Suppose $W$ is a linear $\mathbb{A}$-subspace of $\mathbb{A}^{n}$ and
$\phi:\mathbb{A}^{n}\rightarrow\mathbb{A}^{n}$ is a twisted isomorphism of
$\mathbb{A}^{n}$. Then the image of $\phi$ is also a linear $\mathbb{A}%
$-subspace of $\mathbb{A}^{n}$
\end{lemma}

Proof: This is obvious when $\mathbb{A}$ is either $\mathbb{R}$ or
$\mathbb{C}$ because twisted isomorphisms are the same as vector space
isomorphisms. When $\mathbb{A}=\mathbb{H}$, we have $\phi=\left(  \alpha
,\beta\right)  \in G_{\mathbb{H}}\left(  n\right)  =Sp\left(  n\right)
Sp\left(  1\right)  $. If $\beta=1$ then $\phi$ is an automorphism of
$\mathbb{A}^{n}$, thus it transforms linear $\mathbb{H}$-subspaces to one
another. In any event, the action of $\beta\in Sp\left(  1\right)  $ on
$\mathbb{A}^{n}$ is the diagonal action on the right, thus its also stabilize
any linear $\mathbb{H}$-subspace of $\mathbb{A}^{n}$. Hence we have the
result. $\blacksquare$

\bigskip

From the above lemma, we have the following well-defined notion.

\begin{definition}
Let $M$ be any Riemannian $\mathbb{A}$-manifold. A submanifold $C$ of $M$ is
called a $\mathbb{A}$-submanifold if for any point $p$ in $C$, its tangent
space $T_{p}C$ is a linear $\mathbb{A}$-subspace of $T_{p}M$.
\end{definition}

It is easy to see that any $\mathbb{A}$-submanifold is itself a Riemannian
$\mathbb{A}$-manifold. In the real case, a $\mathbb{R}$-submanifold is simply
an ordinary submanifold. In the complex case, a $\mathbb{C}$-submanifold in a
K\"{a}hler manifold $M$ is equivalent to a complex submanifold of $M$. It
always have absolute minimal volume by the Wirtinger formula, or via
calibration theory.

In the quaternionic case, Gray \cite{Gray} shows that a $\mathbb{H}%
$-submanifold in a quaternionic K\"{a}hler manifold $M$ is always a totally
geodesic submanifold. In particular they are rather rare. Note that a
submanifold $C$ in a hyperk\"{a}hler manifold $M$ which is complex with
respect to all $I$, $J$ and $K$, i.e. a hyperholomorphic submanifold, is a
$\mathbb{H}$-submanifold of $M$.

Remark: If $f$ is an $\mathbb{A}$-isometry of a Riemannian $\mathbb{A}%
$-manifold $M$, then its fixed point set is always an $\mathbb{A}$-submanifold
of $M$.

\subsection{\label{1Sec B Lag submfd}$\frac{1}{2}\mathbb{A}$-Lagrangian submanifolds}

Symplectic geometry is a subject about Lagrangian submanifolds in a symplectic
manifold $M$. If $M$ has a compatible Riemannian metric such that the
symplectic form is parallel, then $M$ is a K\"{a}hler manifold. To rigidity
Lagrangian submanifolds using the metric, we need $M$ to be a Calabi-Yau
manifold and we consider those Lagrangian submanifolds $C$ which are
\textit{special} in the sense that $\operatorname{Im}\Omega|_{C}=0$. They are
an essential ingredient in the mirror symmetry conjecture as discovered by
Strominger, Yau and Zaslow \cite{SYZ}. Recently we also realized that complex
Lagrangian submanifolds in hyperk\"{a}hler manifolds, Cayley submanifold in
$Spin\left(  7\right)  $-manifolds, associative and coassociative submanifolds
in $G_{2}$-manifolds all play important roles in conformal field theory,
string theory and M-theory, and therefore on the geometry of these manifolds.
As we will explain below, all these are (special) $\frac{1}{2}\mathbb{A}%
$-Lagrangian submanifolds and play the role of decomplexifications of
(special) $\mathbb{A}$-manifolds.

\bigskip

First we study the linear case and we begin by reviewing linear Lagrangian
subspaces in $\mathbb{C}^{n}$ with the standard Hermitian complex structure
$J$ and the standard symplectic structure $\omega=\Sigma dx^{j}\wedge dy^{j}$.
A real linear subspace $C$ of dimension $n$ in $\mathbb{C}^{n}$ is called
\textit{Lagrangian} if
\[
\omega|_{C}=0\text{.}%
\]
Equivalently it satisfies
\[
JC\,\text{ is perpendicular to }C\text{,}%
\]
because%
\[
\omega\left(  u,v\right)  =g\left(  Ju,v\right)  \text{.}%
\]
That is we have an orthogonal decomposition
\[
\mathbb{C}^{n}=C\oplus JC
\]
and therefore we can regard a Lagrangian subspace as giving a \textit{real
structure} on $\mathbb{C}^{n}$.

For example $\mathbb{R}^{n}$, the fixed point set of the complex conjugation,
is a Lagrangian subspace in $\mathbb{C}^{n}$. In fact every Lagrangian
submanifold in $\mathbb{C}^{n}$ can be brought to $\mathbb{R}^{n}$ by some
unitary transformation of $\mathbb{C}^{n}$ and the set of all Lagrangian
subspaces in $\mathbb{C}^{n}$ is a homogeneous space $U\left(  n\right)
/SO\left(  n\right)  $.

We are going to generalize this to other normed linear $\mathbb{A}$-spaces
$V\cong\mathbb{A}^{n}$. Note that every imaginary element $u$ in $\mathbb{A}$
with unit length, i.e. $u\in S\left(  \operatorname{Im}\mathbb{A}\right)  $,
defines a complex structure $J_{u}$ on $V$,
\[
J_{u}\left(  y\right)  =yu,
\]
for any $y\in V$. This is because every normed algebra $\mathbb{A}$ is
alternative, $\left(  yu\right)  u=y\left(  u^{2}\right)  $,$\,$and this
implies that $\left(  J_{u}\right)  ^{2}=-1$, i.e. a complex structure on $V$.
Now we define a $\frac{1}{2}\mathbb{A}$-Lagrangian linear subspace in $V$ to
be a\textit{\ maximally real }subspace in $V$.

\begin{definition}
Suppose $C$ is a middle dimensional real linear subspace in a normed linear
$\mathbb{A}$-space $V$. It is called a $\frac{1}{2}\mathbb{A}$-Lagrangian
linear subspace in $V$ if
\[
J_{u}C\text{ }\bot\text{ }C
\]
for any unit element $u\in L$ in some real linear subspace $L\subset
\operatorname{Im}\mathbb{A}$ of $\dim L=\frac{1}{2}\dim\mathbb{A}$.
\end{definition}

This definition is equivalent to having an orthogonal decomposition
\[
V=C\oplus J_{u}C,
\]
for every $u\in L$ with $\left|  u\right|  =1$.

\begin{theorem}
If $C$ is a $\frac{1}{2}\mathbb{A}$-Lagrangian linear subspace in $V$ and we
write $L\subset\operatorname{Im}\mathbb{A}$ as in the above definition, then
$C$ is a $J_{v}$-complex linear subspace of $V$ for any unit vector
$v\in\operatorname{Im}\mathbb{A}$ perpendicular to $L$.
\end{theorem}

Proof: This is an empty assertion when $\mathbb{A}$ is either $\mathbb{R}$ or
$\mathbb{C}$. When $\mathbb{A}=\mathbb{H}$ with the standard complex
structures $I$, $J$, $K$, we can assume, without loss of generality,
\[
IC\,\bot\,\,C\text{ and }JC\,\,\bot\,\,C\text{.}%
\]
Using the associativity of $\mathbb{H}$ and $C$ being of middle dimensional,
we have
\[
KC=\left(  IJ\right)  C=I\left(  JC\right)  =C.
\]
That is $C$ is a complex $K$-linear subspace of $V$, thus proving the assertion.

When $\mathbb{A}=\mathbb{O}$, the above arguments do not work since
$\mathbb{O}$ is non-associative. First we claim that we can assume that the
$\frac{1}{2}\mathbb{A}$-Lagrangian linear subspace $C$ in $V\cong\mathbb{O}$
contains $\mathbb{R}$. To see this, suppose that $u\in\operatorname{Im}%
\mathbb{O}$ with $\left|  u\right|  =1$ satisfies $J_{u}C\,\bot C$, i.e.
\[
\left\langle cu,c^{\prime}\right\rangle =0
\]
for any $c,c^{\prime}\in C$. We recall that for any $g\in Spin\left(
7\right)  $ there is an $\theta\in SO\left(  8\right)  $ such that $g\left(
cu\right)  =g\left(  c\right)  \theta\left(  u\right)  $. In particular,
$u\in\operatorname{Im}\mathbb{O}$ implies that $\theta\left(  u\right)  $ is
also imaginary because
\[
0=\left\langle u,1\right\rangle =\left\langle g\left(  1\right)  \theta\left(
u\right)  ,g\left(  1\right)  \right\rangle =\left|  g\left(  1\right)
\right|  ^{2}\left\langle \theta\left(  u\right)  ,1\right\rangle \text{.}%
\]

From $\left\langle cu,c^{\prime}\right\rangle =0$, we obtain
\[
\left\langle g\left(  c\right)  \theta\left(  u\right)  ,g\left(  c^{\prime
}\right)  \right\rangle =0,
\]
with $\theta\left(  u\right)  \in\operatorname{Im}\mathbb{O}$, $\left|
\theta\left(  u\right)  \right|  =1$. That is
\[
J_{\theta\left(  u\right)  }g\left(  C\right)  \,\bot\,g\left(  C\right)
\text{.}%
\]
Using the fact that $Spin\left(  7\right)  $ acts transitively on the unit
sphere in $\mathbb{O}$, we have therefore verified our claim.

Since $\mathbb{R}\subset C$, we have an orthogonal direct sum decomposition,
\[
C=\mathbb{R}\oplus D
\]
for some three dimension subspace $D\subset\operatorname{Im}\mathbb{O}$.

Second we claim that there is an orthogonal direct sum decomposition
\[
\operatorname{Im}\mathbb{O}=D\oplus L\text{.}%
\]
The reason is, for any $d\in D$ and $u\in L$ with $\left|  u\right|  =1$, we
have
\[
\left\langle u,d\right\rangle =\left\langle 1,d\bar{u}\right\rangle
=-\left\langle 1,J_{u}\left(  d\right)  \right\rangle =0
\]
because $1\in C$ which is perpendicular to $J_{u}\left(  d\right)  \in J_{u}C$.

Third we want to prove that $C$ is the real span of $1$, $i$, $j$ and $k$,
possibly after a $G_{2}$ rotation. Without loss of generality, we can assume
that $i,j\in D$, by a $G_{2}$ rotation if necessary. We write
\[
k=ij=d+u\in\operatorname{Im}\mathbb{O}%
\]
for some element $d\in D$ perpendicular to both $i$ and $j$, and for some
$u\in L$. Then
\[
\left\langle u,u\right\rangle =\left\langle u,k\right\rangle -\left\langle
u,d\right\rangle =\left\langle u,ij\right\rangle \text{,}%
\]
because $d$ and $u$ are perpendicular by the second claim. On the other hand,
\[
\left\langle u,ij\right\rangle =\left\langle \bar{\imath}u,j\right\rangle
=\left\langle \bar{\imath},j\bar{u}\right\rangle =\left\langle
i,ju\right\rangle =0,
\]
because $i\in C$ and $ju\in Cu$ are perpendicular to each other. This implies
that $u=0$, or equivalently, $k\in D$. That is $C$ is the real span of $1$,
$i$, $j$, $k$ and $L$ is the real span of $e$, $ei$, $ej$, $ek$ with
$\mathbb{O}=\mathbb{H}\oplus e\mathbb{H}$. This, in particular, gives our
theorem. $\blacksquare$

\bigskip

Remark: The converse to the above corollary is not true. For example
$C=\mathbb{H}\times\left\{  0\right\}  \subset\mathbb{H}^{2}$ is $J$-linear
but it is not a $\frac{1}{2}\mathbb{A}$-Lagrangian subspace of $V=\mathbb{H}%
^{2}$.

From the proof, we see that a $\mathbb{H}$-Lagrangian linear subspace in
$\mathbb{O}$ is precisely $\mathbb{H\times}\left\{  0\right\}  \subset
\mathbb{H}^{2}=\mathbb{O}$ up to the action of $Spin\left(  7\right)  $. Such
linear subspaces are studied by Harvey and Lawson and they are called
\textit{Cayley subspaces} of $\mathbb{O}$ and they are volume calibrated by
the four form $\Theta$ on $\mathbb{O}$ (see \cite{HL} for various
characterizations of Cayley subspaces). We therefore have the following corollary.

\begin{corollary}
$\mathbb{H}$-Lagrangian linear subspaces in $\mathbb{O}$ are equivalent to
Cayley subspaces.
\end{corollary}

Moreover we can justify the definition of $\frac{1}{2}\mathbb{A}$-Lagrangian
linear subspace as being maximally real subspaces of $V$, in the following lemma.

\begin{corollary}
Suppose that $C$ is a middle dimensional real linear subspace in a normed
linear $\mathbb{A}$-space $V$. If there is a real linear subspace $L^{\prime
}\subset\operatorname{Im}\mathbb{A}$ such that
\[
J_{u}C\text{ }\bot\text{ }C
\]
for any unit element $u\in L^{\prime}$, then $\dim L^{\prime}\leq\frac{1}%
{2}\dim\mathbb{A}$.
\end{corollary}

\bigskip

For submanifolds in a K\"{a}hler manifold, the condition $JC$ $\bot$ $C$ is
usually written as $\omega|_{C}=0$, namely the Lagrangian condition. We want
to do the same thing for all $\frac{1}{2}\mathbb{A}$-Lagrangian linear subspaces.

\begin{proposition}
There is an embedding of $\operatorname{Im}\mathbb{A}$ into the space of two
forms on the real vector space $V_{\mathbb{R}}$,
\begin{gather*}
\operatorname{Im}\mathbb{A}\overset{\subset}{\rightarrow}\Lambda
^{2}V_{\mathbb{R}}^{\ast},\\
u\rightarrow\omega_{u}%
\end{gather*}
defined by
\[
\omega_{u}\left(  x\otimes y\right)  =\left\langle x,yu\right\rangle \text{.}%
\]
\end{proposition}

\begin{lemma}
Moreover this embedding intertwine the action of $G_{\mathbb{A}}\left(
n\right)  $ on $\Lambda^{2}V_{\mathbb{R}}^{\ast}$ via the inclusion
$G_{\mathbb{A}}\left(  n\right)  \subset O\left(  m\right)  ,$ and a natural
action of $G_{\mathbb{A}}\left(  n\right)  $ on $\operatorname{Im}\mathbb{A}$
given as follows: When $\mathbb{A}$ $=\mathbb{C\,}$or $\mathbb{H}$ the action
is given by the composition of $\lambda_{\mathbb{A}}\left(  n\right)  $ and
the conjugation; when $\mathbb{A}=\mathbb{O}$, the action is given by
$Spin\left(  7\right)  \rightarrow SO\left(  \operatorname{Im}\mathbb{O}%
\right)  $, $g\rightarrow\theta_{g}$.
\end{lemma}

Proof: First $\omega_{u}$ is a two form because when $u\in\operatorname{Im}%
\mathbb{O}$, we have%

\[
\omega_{u}\left(  x\otimes x\right)  =\left\langle x,xu\right\rangle =\left|
x\right|  ^{2}\left\langle 1,u\right\rangle =0\text{.}%
\]

Next we prove the compatibility with respect to the actions by $G_{\mathbb{A}%
}\left(  n\right)  $. In the complex case, it follows from $G_{\mathbb{A}%
}\left(  n\right)  =U\left(  n\right)  $ acts trivially on $\operatorname{Im}%
\mathbb{C}$. In the quaternionic case, $\left(  \alpha,\beta\right)  \in
Sp\left(  n\right)  Sp\left(  1\right)  $ acts on $u\in\operatorname{Im}%
\mathbb{H}$ and gives $\beta^{-1}u\beta$. We compute,%
\[
\left\langle \alpha x\beta,\alpha y\beta\left(  \beta^{-1}u\beta\right)
\right\rangle =\left\langle x\beta,yu\beta\right\rangle =\left\langle
x,yu\right\rangle ,
\]
and we have the claim in this case. For the octonionic case, for any $g\in
Spin\left(  7\right)  $, we have%
\[
g\left(  yu\right)  =g\left(  y\right)  \theta_{g}\left(  u\right)
\]
for any $y,u\in\mathbb{O}$. We compute%
\[
\omega_{u}\left(  x\otimes y\right)  =\left\langle x,yu\right\rangle
=\left\langle g\left(  x\right)  ,g\left(  yu\right)  \right\rangle
=\left\langle g\left(  x\right)  ,g\left(  y\right)  \theta_{g}\left(
u\right)  \right\rangle =\omega_{\theta_{g}\left(  u\right)  }\left(  g\left(
x\right)  \otimes g\left(  y\right)  \right)  \text{,}%
\]
and hence the result. $\blacksquare$

\bigskip

Remark: This proposition implies that, on any Riemannian $\mathbb{A}$-manifold
$M$, there is a real vector subbundle over $\Lambda^{2}\left(  M\right)  $ of
rank $\dim\mathbb{A}-1$, denotes $\mathbf{s}_{\mathbb{A}}\left(  T_{M}\right)
$, such that each fiber is a copy of $\operatorname{Im}\mathbb{A}$.

\bigskip

We can identify $\operatorname{Im}\mathbb{A}\subset\Lambda^{2}V_{\mathbb{R}%
}^{\ast}\cong\mathbf{o}\left(  m\right)  $ explicitly. Obviously
$\operatorname{Im}\mathbb{R}=0$. In the complex case, if we choose $u=i$ then
$\omega_{u}$ is simply the standard K\"{a}hler form as can be easily checked.
Similarly, in the quaternionic case, $\omega_{i}$, $\omega_{j}$ and
$\omega_{k}$ are the standard K\"{a}hler forms for the complex structures $I$,
$J$ and $K$ respectively. In particular we have the following decomposition,
\[
\mathbf{g}_{\mathbb{A}}\left(  n\right)  =\mathbf{h}_{\mathbb{A}}\left(
n\right)  +\operatorname{Im}\mathbb{A},
\]
for $\mathbb{A}\neq\mathbb{O}$. In the octonionic case, we have
\[
\operatorname{Im}\mathbb{A}\cong\Lambda_{7}^{2}\mathbb{O}%
\]
where
\[
\Lambda^{2}\mathbb{O}=\Lambda_{21}^{2}\mathbb{O}+\Lambda_{7}^{2}%
\mathbb{O}\text{,}%
\]
is the decomposition of $\Lambda^{2}\mathbb{O}$ into irreducible $Spin\left(
7\right)  $-representations. This can be verified either by identifying the
$Spin\left(  7\right)  $representation $\Lambda_{7}^{2}\mathbb{O}$ as given by
$g\rightarrow\theta_{g}$, or simply by checking the dimensions of these
irreducible pieces. Thus we have obtained the first part of the following proposition.

\begin{proposition}
Suppose that $V$ is a normed linear $\mathbb{A}$-space and we denote the image
of $\operatorname{Im}\mathbb{A}$ in $\Lambda^{2}V_{\mathbb{R}}^{\ast}$ as
$\mathbf{s}_{\mathbb{A}}$. There is a decomposition of $\mathbf{g}%
_{\mathbb{A}}\left(  n\right)  $-representations,
\[
\mathbf{g}_{\mathbb{A}}\left(  n\right)  =\mathbf{h}_{\mathbb{A}}\left(
n\right)  +\mathbf{s}_{\mathbb{A}},
\]
when $\mathbb{A\neq O}$ and $\Lambda^{2}\mathbb{O}=\Lambda_{21}^{2}%
\mathbb{O}+\mathbf{s}_{\mathbb{O}}$.

Moreover for any $u\in\operatorname{Im}\mathbb{A}$ with unit length, we have a
natural complex structure $J_{u}$ on $V$ and $\omega_{u}\in\mathbf{s}%
_{\mathbb{A}}\subset\Lambda^{2}V_{\mathbb{R}}^{\ast}$ satisfying
\[
V=C\oplus J_{u}C\text{ if and only if }\omega_{u}|_{C}=0.
\]
\end{proposition}

Proof: The first half is proven above. The proof for the second half is the
same as in the K\"{a}hler case, namely for any $x,y\in C$, we have
\[
\omega_{u}\left(  x\wedge y\right)  =\left\langle x,yu\right\rangle
=\left\langle x,J_{u}y\right\rangle .
\]
The claim follows from the middle dimensionality of $C$. $\blacksquare$

\bigskip

As a corollary, we obtain an equivalent definition of a $\frac{1}{2}%
\mathbb{A}$-Lagrangian linear subspace which resemblance the definition of an
ordinary Lagrangian subspace.

\begin{corollary}
Suppose that $C$ is a half dimensional real linear subspace in a normed linear
$\mathbb{A}$-space $V$. Then $C$ is a $\frac{1}{2}\mathbb{A}$-Lagrangian
linear subspace if and only if there is a linear subspace $L\subset
\operatorname{Im}\mathbb{A}\cong\mathbf{s}_{\mathbb{A}}\subset\Lambda
^{2}V_{\mathbb{R}}^{\ast}$ with $\dim L=\frac{1}{2}\dim\mathbb{A}$ such that%
\[
\omega|_{C}=0
\]
for any $\omega\in L$.
\end{corollary}

\bigskip

Using the above proposition and the proof of the earlier theorem which
describe explicitly $\mathbb{H}$-Lagrangian linear subspaces in $\mathbb{O}$,
we obtain the following corollary.

\begin{corollary}
If $C$ is a real four dimensional linear subspace of $\mathbb{O}$, then it is
a $\mathbb{H}$-Lagrangian linear subspace if and only if the homomorphism
defined by restricting differential forms,
\[
\Lambda_{7}^{2}\left(  \mathbb{O}\right)  \rightarrow\Lambda^{2}\left(
C\right)
\]
has a four dimensional kernel. Moreover this happens exactly when the image of
the above homomorphism is $\Lambda_{+}^{2}\left(  C\right)  $.
\end{corollary}

The next proposition is basically well-known (see \cite{HL} for the proof in
the octonion case).

\begin{proposition}
The space of $\frac{1}{2}\mathbb{A}$-Lagrangian linear subspaces in
$\mathbb{A}^{n}$ is a homogeneous space of $G_{\mathbb{A}}\left(  n\right)  $.
Explicitly they are $U\left(  n\right)  /SO\left(  n\right)  $, $Sp\left(
n\right)  Sp\left(  1\right)  /U\left(  n\right)  U\left(  1\right)  $ and
$Spin\left(  7\right)  /Sp\left(  1\right)  ^{3}$ for $\mathbb{A}=\mathbb{C}$,
$\mathbb{H}$ and $\mathbb{O}$ respectively.
\end{proposition}

\bigskip

After all these studies of $\frac{1}{2}\mathbb{A}$-Lagrangian in the linear
case, we come to the definition of a $\frac{1}{2}\mathbb{A}$-Lagrangian submanifold.

\begin{definition}
A middle dimensional real submanifold $C$ in a Riemannian $\mathbb{A}%
$-manifold $M$ is called a $\frac{1}{2}\mathbb{A}$-Lagrangian submanifold of
$M$ if there is a vector subbundle\footnote{In fact we only need the subbundle
$L$ to be defined over $C$.} $L\subset\mathbf{s}_{\mathbb{A}}\left(
T_{M}\right)  \subset\Lambda^{2}\left(  M\right)  $ of $rank\left(  L\right)
=\frac{1}{2}\dim\mathbb{A}$ such that
\[
\omega|_{C}=0
\]
for any smooth section $\omega\in\Gamma\left(  M,L\right)  $.
\end{definition}

\bigskip

Remark: Unlike Lagrangian submanifolds in a K\"{a}hler manifold, $\mathbb{C}%
$-Lagrangian submanifolds in a quaternionic K\"{a}hler manifold are not widely
studied in the literature. Nonetheless, every surface in an oriented four
manifold is a $\mathbb{C}$\--Lagrangian submanifold.

The following table summarizes $\frac{1}{2}\mathbb{A}$-Lagrangian submanifolds
in Riemannian $\mathbb{A}$-manifolds and their common names.%

\[%
\begin{tabular}
[c]{|l|l|}\hline
$%
\begin{array}
[c]{l}%
\,\\
\,
\end{array}
$ & $\frac{1}{2}\mathbb{A}$-Lagrangian submanifolds\\\hline\hline
$%
\begin{array}
[c]{l}%
\,\\
\,
\end{array}
\mathbb{C}$ & $\omega|_{C}=0$\\
& (Lagrangian submanifolds)\\\hline
$%
\begin{array}
[c]{l}%
\,\\
\,
\end{array}
\mathbb{H}$ & $\omega|_{C}=0$ for all \textbf{\ }$\omega\in L\subset Sym^{2}S$
w/ rank$\left(  L\right)  =2$\\
& \\\hline
$%
\begin{array}
[c]{l}%
\,\\
\,
\end{array}
\mathbb{O}$ & $\omega|_{C}=0$ for all \textbf{\ }$\omega\in L\subset
\Lambda_{7}^{2}\left(  M\right)  $ w/ rank$\left(  L\right)  =4$\\
& (Cayley submanifolds)\\\hline
\end{tabular}
\]

\bigskip

\subsection{\label{1Sec Sp B Lag submfd}Special $\frac{1}{2}\mathbb{A}%
$-Lagrangian submanifolds}

Next we introduce \textit{special} $\frac{1}{2}\mathbb{A}$-Lagrangian
submanifolds inside special $\mathbb{A}$-manifolds. As before, we start with
the linear theory about special $\frac{1}{2}\mathbb{A}$-Lagrangian linear
subspaces in $\mathbb{A}^{n}$.

For example in the complex case, if $C$ is a Lagrangian linear subspace of
$V\cong\mathbb{C}^{n}$, then there is an $A\in SU\left(  n\right)  $ such
that
\[
A\left(  C\right)  =e^{i\theta}\mathbb{R}^{n}%
\]
for some angle $\theta$, which is usually called the \textit{phase} of $C$.
When $C$ is a Lagrangian submanifold in a Calabi-Yau manifold $M$, then the
gradient of $\theta$ is the mean curvature vector of $C$ in $M$. In fact, if
$\theta$ is constant over $C$, then $C$ is a minimal submanifold in $M$ with
absolute minimal volume and it is called a special Lagrangian submanifold in
$M$ with phase $\theta$. We want to generalize this concept to other special
$\mathbb{A}$-manifolds. First we need the following definitions, which are
roughly the \textit{determinants} of $C$.

\begin{definition}
For any $\frac{1}{2}\mathbb{A}$-Lagrangian subspace $C$ in $V\cong
\mathbb{A}^{n}$ we define another $\frac{1}{2}\mathbb{A}$-Lagrangian subspace
$\lambda\left(  C\right)  $ in $\mathbb{A}$ as follows:

When $\mathbb{A}=\mathbb{C}$, $\Lambda\left(  C\right)  $ is the image of
$\Lambda^{n}C$ under the following composition of natural homomorphisms,
\[
\Lambda^{n}V_{\mathbb{R}}\rightarrow\left(  \Lambda^{n,0}\left(  V\right)
\oplus\Lambda^{0,n}\left(  V\right)  \right)  _{\mathbb{R}}\cong\Lambda
^{n,0}\left(  V\right)  =\Lambda^{n}\left(  V\right)  \cong\mathbb{C}\text{.}%
\]
Here the first homomorphism is the orthogonal projection to the Hodge $\left(
p,q\right)  $-decomposition.

When $\mathbb{A}=\mathbb{H}$, any $\mathbb{C}$-Lagrangian submanifold is
holomorphic with respect to a unique complex structures $\pm J_{u}$ with
$u\in\operatorname{Im}\mathbb{H}$. $\lambda\left(  C\right)  $ is the complex
line in $\mathbb{H\cong C}^{2}$ corresponding to $u$ under the following
identification,
\[
S\left(  \operatorname{Im}\mathbb{H}\right)  /\mathbb{Z}_{2}\cong S\left(
\Lambda_{+}^{2}\mathbb{H}\right)  /\mathbb{Z}_{2}\cong\mathbb{P}_{\mathbb{C}%
}\left(  \mathbb{C}^{2}\right)  /\mathbb{Z}_{2}\text{.}%
\]
$\lambda\left(  C\right)  \subset\mathbb{H}$ is only well-defined up to
replacing it by its orthogonal complement $\lambda\left(  C\right)  ^{\bot
}\subset\mathbb{H}$.

When $\mathbb{A}=\mathbb{O}$, we simply define $\lambda\left(  C\right)  =C$.
\end{definition}

Remark: In the complex case, $C$ being a Lagrangian in $V$ implies that
$\lambda\left(  C\right)  $ is a line in $\Lambda^{n}V\cong\mathbb{C}$, thus
also a $\mathbb{R}$-Lagrangian by trivial reason. It can be checked directly
that this is the line in $\mathbb{C}$ with slope $\tan\theta$.

\bigskip

\begin{definition}
A $\frac{1}{2}\mathbb{A}$-Lagrangian linear subspace $C$ in a normed linear
$\mathbb{A}$-space $V\cong\mathbb{C}^{n}$ is called special of type I (resp.
type II) if $1\in\lambda\left(  C\right)  \subset\mathbb{A}$ (resp.
$1\in\lambda\left(  C\right)  ^{\bot}\subset\mathbb{A}$).
\end{definition}

Remark: When $\mathbb{A}=\mathbb{H}$, type I and type II $\mathbb{C}%
$-Lagrangian subspaces are equivalent due to quotient by $\left\{
\pm1\right\}  $ in $G_{\mathbb{H}}\left(  n\right)  =Sp\left(  n\right)
\times Sp\left(  1\right)  /\pm1$.

Remark: As in the $\frac{1}{2}\mathbb{A}$-Lagrangian case, the space of
special $\frac{1}{2}\mathbb{A}$-Lagrangian linear subspaces in $\mathbb{A}%
^{n}$ is a homogeneous space of $H_{\mathbb{A}}\left(  n\right)  $. Explicitly
they are the fibers of the following fiber bundles (see \cite{HL} for the
proof in the octonion case),
\[%
\begin{array}
[c]{llllll}%
\frac{SU\left(  n\right)  }{SO\left(  n\right)  } & \rightarrow &
\frac{U\left(  n\right)  }{SO\left(  n\right)  } & \rightarrow &
\frac{U\left(  1\right)  }{SO\left(  1\right)  } & =S^{1}\\
&  &  &  &  & \\
\frac{Sp\left(  n\right)  }{U\left(  n\right)  } & \rightarrow &
\frac{Sp\left(  n\right)  Sp\left(  1\right)  }{U\left(  n\right)  U\left(
1\right)  }\  & \rightarrow & \frac{Sp\left(  1\right)  }{U\left(  1\right)
}/\pm1 & =S^{2}/\pm1\\
&  &  &  &  & \\
\frac{G_{2}}{Sp\left(  1\right)  ^{2}} & \rightarrow & \frac{Spin\left(
7\right)  }{Sp\left(  1\right)  ^{3}} & \rightarrow & \frac{S^{7}}{Sp\left(
1\right)  } & =S^{4}.
\end{array}
\]

\bigskip

To define the corresponding notion for submanifolds, we recall that there is a
canonical $\mathbb{A}$-bundle $\mathbb{A}_{M}$ over any Riemannian
$\mathbb{A}$-manifold $M$ corresponding to the representation $\lambda
_{\mathbb{A}}\left(  n\right)  :G_{\mathbb{A}}\left(  n\right)  \rightarrow
O\left(  \mathbb{A}\right)  $, which is trivial when $M$ is special. We fix a
trivialization compatible with the action by $\mathbb{A}$ and let $s$ be the
section of this bundle corresponding to $1\in\mathbb{A}$. From the above
linear considerations, if $C$ is a $\frac{1}{2}\mathbb{A}$-Lagrangian
submanifold in $M$, the $\lambda\left(  C\right)  $ is a subbundle of
$\mathbb{A}_{M}$ restricted to $C$.

\begin{definition}
A $\frac{1}{2}\mathbb{A}$-Lagrangian submanifold $C$ in a special Riemannian
$\mathbb{A}$-manifold $M$ is called special of type I (resp. type II) if
$s\in\lambda\left(  C\right)  \subset\mathbb{A}_{M}$ (resp. $s\in
\lambda\left(  C\right)  ^{\bot}\subset\mathbb{A}_{M}$).
\end{definition}

Note that $\mathbb{C}$-Lagrangians of type I and type II are the same.

We have the following characterizations of special $\frac{1}{2}\mathbb{A}%
$-Lagrangian submanifolds, together with their common names.%

\[%
\begin{tabular}
[c]{|l|l|l|}\hline
$%
\begin{array}
[c]{l}%
\,\\
\,
\end{array}
$ & Special $\frac{1}{2}\mathbb{A}$-Lagr. submfd. (type I) & Special $\frac
{1}{2}\mathbb{A}$-Lagr. submfd. (type II)\\\hline\hline
$\mathbb{C}$ & $\omega|_{C}=\operatorname{Im}\Omega|_{C}=0
\begin{array}
[c]{l}%
\,\\
\,
\end{array}
$ & $\omega|_{C}=\operatorname{Re}\Omega|_{C}=0$\\
& (special Lagr. submfd. with phase $0$) & (special Lagr. submfd. with phase
$\pi/2$)\\\hline
$\mathbb{H}$ & $\Omega|_{C}=0
\begin{array}
[c]{l}%
\,\\
\,
\end{array}
$ & $\Omega|_{C}=0
\begin{array}
[c]{l}%
\,\\
\,
\end{array}
$\\
& (Complex Lagrangian submanifolds) & (Complex Lagrangian
submanifolds)\\\hline
$\mathbb{O}$ & $\times$ preserves $C$ (or $\chi|_{C}=0$)$%
\begin{array}
[c]{l}%
\,\\
\,
\end{array}
$ & $\Omega|_{C}=0$\\
& (Associative submanifolds) & (Coassociative submanifolds)\\\hline
\end{tabular}
\]
\ 

\subsection{\label{1Sec Vol calib}Volume calibrations}

This short subsection is included for completeness, readers should consult the
paper \cite{HL} by Harvey and Lawson for a careful treatment.

\begin{definition}
(1) A differential form $\Phi\in\Omega^{k}\left(  M\right)  $ in an oriented
Riemannian manifold $M$ is called a volume calibrating form if
\[%
\begin{array}
[c]{cc}%
d\Phi=0 & \\
\Phi|_{P}\leq dv_{P} & \text{ for all }P\in\widetilde{Gr}\left(
k,T_{M}\right)  \text{.}%
\end{array}
\]

(2) A $k$-dimensional submanifold $C\subset M$ is calibrated by $\Phi$ if
\[
\Phi|_{C}=dv_{C}\text{.}%
\]
\end{definition}

We have the following fundamental lemma.

\begin{lemma}
Any closed calibrated submanifold $C$ is homologically volume minimizing, i.e.
$Vol\left(  C\right)  \leq Vol\left(  C^{\prime}\right)  $ provided $C$ and
$C^{\prime}$ represent the same homology class in $M$. Moreover if $Vol\left(
C\right)  =Vol\left(  C^{\prime}\right)  $ then $C^{\prime}$ is also calibrated.
\end{lemma}

For Riemannian $\mathbb{A}$-manifolds, there are natural calibrating forms
$\Phi$ and their calibrating submanifolds are closely related to $\mathbb{A}%
$-submanifolds. We list them in the following table.%

\[%
\begin{tabular}
[c]{|l|l|l|}\hline
$%
\begin{array}
[c]{l}%
\,\\
\,
\end{array}
$ & Holonomy & Calibrating form\\
& (Riemannian $\mathbb{A}$-manifolds) & (Calibrated
submanifolds)\\\hline\hline
$\mathbb{R}$ & $O\left(  n\right)
\begin{array}
[c]{l}%
\,\\
\,
\end{array}
$ & $\Phi=\exp\left(  \nu ol_{M}\right)  $\\
& (Riemannian manifolds) & (Points and $M$)\\\hline
$\mathbb{C}$ & $U\left(  n\right)
\begin{array}
[c]{l}%
\,\\
\,
\end{array}
$ & $\Phi=\exp\left(  \omega\right)  $\\
& (Kahler manifolds)\thinspace & (Complex submanifolds)\\\hline
$\mathbb{H}$ & $Sp\left(  n\right)  Sp\left(  1\right)
\begin{array}
[c]{l}%
\,\\
\,
\end{array}
$ & $\Phi=\exp\left(  \Theta\right)  $\\
& (Quaternionic Kahler mfds) & (Quaternionic submanifolds)\\\hline
$\mathbb{O}$ & $Spin\left(  7\right)  $ & $\Phi=\exp\left(  \Theta\right)
\begin{array}
[c]{l}%
\,\\
\,
\end{array}
$\\
& ($Spin\left(  7\right)  $-manifolds) & (Cayley submanifods)\\\hline
\end{tabular}
\]

For special Riemannian $\mathbb{A}$-submanifolds, there are further
calibrating forms, which are closely related to special $\frac{1}{2}%
\mathbb{A}$-Lagrangian submanifolds above. They are listed in the next table.
\[%
\begin{tabular}
[c]{|l|l|l|}\hline
$%
\begin{array}
[c]{l}%
\,\\
\,
\end{array}
$ & Holonomy & Calibrating form\\
& (Special Riem. $\mathbb{A}$-manifolds) & (Calibrated
submanifolds)\\\hline\hline
$\mathbb{R}$ & $SO\left(  n\right)
\begin{array}
[c]{l}%
\,\\
\,
\end{array}
$ & $\Phi=v_{M}$\\
& (Oriented manifolds) & (Whole manifold $M$)\\\hline
$\mathbb{C}$ & $SU\left(  n\right)
\begin{array}
[c]{l}%
\,\\
\,
\end{array}
$ & Type I: $\Phi=\operatorname{Re}\Omega$\\
& (Calabi-Yau manifolds)\thinspace & (SLag with phase $0$)\\
& $%
\begin{array}
[c]{l}%
\,\\
\,
\end{array}
$ & Type II: $\Phi=\operatorname{Im}\Omega$\\
&  & (SLag with phase $\pi/2$)\\\hline
$\mathbb{H}$ & $Sp\left(  n\right)
\begin{array}
[c]{l}%
\,\\
\,
\end{array}
$ & $\Phi_{1}=\operatorname{Re}\Omega_{I}^{n}$ and $\Phi_{2}=\operatorname{Re}%
\Omega_{K}^{n}$\\
& (Hyperkahler manifolds) & (Complex Lagrangian submanifolds)\\\hline
$\mathbb{O}$ & $G_{2}$ & Type I: $\Phi=\Omega%
\begin{array}
[c]{l}%
\,\\
\,
\end{array}
$\\
& ($G_{2}$-manifolds) & (Associative submanifods)\\
&  & Type II: $\Phi=\Theta%
\begin{array}
[c]{l}%
\,\\
\,
\end{array}
$\\
&  & (Coassociative submanifolds)\\\hline
\end{tabular}
\]

Remark: A middle dimension submanifold $C$ in a hyperk\"{a}hler manifold $M$
is a complex Lagrangian if $\Omega_{J}|_{C}=0$. Since $\Omega_{J}=\omega
_{I}+i\omega_{K}$ we have $\omega_{I}=\omega_{K}=0.$ In this case $\Omega
_{I}|_{C}=\omega_{J}$ (since $\omega_{K}=0$), so $\operatorname{Im}\Omega
_{I}^{n}|_{C}=0$, i.e. $C$ is calibrated by $\operatorname{Re}\Omega_{I}^{n}$.
Similar for $\operatorname{Re}\Omega_{K}^{n}$. The converse is also true. It
is also calibrated by $\omega_{J}^{n}$.

\section{\label{1Sec Geom Dual}Geometry and duality}

On a $\mathbb{C}$-manifold $M$ (i.e. K\"{a}hler manifold), its $\mathbb{C}%
$-geometry studies cycles $\left(  C,D_{E}\right)  $ with $C$ a complex
submanifold in $M$ and $E$ a holomorphic bundle over $C$. In algebraic
geometry, one also allow $C$ and $E$ to be \textit{singular }and consider
$D^{b}\left(  M\right)  $ the derived category of coherent sheaves on $M$.
When $M$ is special (i.e. Calabi-Yau manifold) we would also require $D_{E}$
to be a special $\mathbb{C}$-connection, i.e. a Hermitian Yang-Mills
connection over $C$.

On the other hand, the $\mathbb{R}$-Lagrangian geometry of a K\"{a}hler (or
symplectic) manifold $M$ studies cycles $\left(  C,D_{E}\right)  $ with $C$ a
Lagrangian submanifold and $D_{E}$ a unitary flat connection over $C$. The
space of morphisms between these cycles are the Floer homology groups. When
$M$ is special, we also study special Lagrangians $C$. For instance the mirror
Calabi-Yau manifold is conjectured to be the moduli space of certain special
Lagrangian cycles, as in the SYZ mirror conjecture.

A novelty about geometry for manifolds with special holonomy is the
\textit{duality} transformation. For example, the mirror symmetry among
Calabi-Yau manifolds (see e.g. \cite{Le Geom MS}), motivated from physics, is
still very mysterious to mathematicians. From the work of Strominger, Yau and
Zaslow \cite{SYZ}, we expect that it is a fiberwise Fourier transformation
along a special Lagrangian torus fibration\footnote{We also need a Legendre
transformation along the base \cite{Le MS corr}.}.

Fourier transformation on tori are well-studied in mathematics. We review it
from our point of view: namely it should transform the $\mathbb{A}^{\prime}%
$-geometry of one torus to the $\mathbb{A}^{\prime}$-geometry of its dual
torus, with $\mathbb{A}^{\prime}=\mathbb{R},\mathbb{C}$ or $\mathbb{H}$. If we
complexify a torus $T^{n}$ to $M=T^{n}\times i\mathbb{R}^{n}$, then it is a
special $\mathbb{A}$-manifold with $\mathbb{A}\,^{\prime}=\frac{1}%
{2}\mathbb{A}$. Moreover it has a natural fibration by special $\frac{1}%
{2}\mathbb{A}$-Lagrangian tori given by projection. The fiberwise Fourier
transformation, or the SYZ transformation, should transform the $\mathbb{A}%
$-geometry on $M$ to the $\frac{1}{2}\mathbb{A}$-Lagrangian geometry on
$W=T^{n\ast}\times i\mathbb{R}^{n\ast}$, and vice versa.

\subsection{\label{1Sec A and A2 geom}$\mathbb{A}$-geometry and $\frac{1}%
{2}\mathbb{A}$-Lagrangian geometry}

Suppose $M$ is a Riemannian $\mathbb{A}$-manifold. We consider (i) the
geometry of $\mathbb{A}$-cycles and (ii) the geometry of $\frac{1}%
{2}\mathbb{A}$-Lagrangian cycles on $M$.

\begin{definition}
Suppose $M$ is a Riemannian $\mathbb{A}$-manifold. A pair $\left(
C,D_{E}\right)  $ is called an (i) $\mathbb{A}$-cycle if $C$ is a $\mathbb{A}%
$-submanifold in $M$ and $D_{E}$ is a $\mathbb{A}$-connection over $C$ or (ii)
$\frac{1}{2}\mathbb{A}$-Lagrangian cycle if $C$ is a $\frac{1}{2}\mathbb{A}%
$-Lagrangian submanifold in $M$ and $D_{E}$ is a special $\frac{1}%
{2}\mathbb{A}$-connection over $C$.
\end{definition}

The following table gives the common names of these cycles.
\[%
\begin{tabular}
[c]{|l||l||l|}\hline
$%
\begin{array}
[c]{l}%
\,\\
\,
\end{array}
$ & $\mathbb{A}$-cycle on $M$ & $\frac{1}{2}\mathbb{A}$-Lagrangian cycle on
$M$\\\hline\hline
$\mathbb{C}$ & Complex submanifold + & Lagrangian submanifold + $%
\begin{array}
[c]{l}%
\,\\
\,
\end{array}
$\\
& Holomorphic bundle & Unitary flat bundle\\\hline
$\mathbb{H}$ & Quaternionic submanifold + $%
\begin{array}
[c]{l}%
\,\\
\,
\end{array}
$ & $\mathbb{C}$-Lagrangian submanifold +\\
& Bundle with $\mathbb{C}$-connection & Holomorphic bundle\\\hline
$\mathbb{O}$ & The whole manifold $M$ + & Cayley submanifold + $%
\begin{array}
[c]{l}%
\,\\
\,
\end{array}
$\\
& $Spin\left(  7\right)  $-Donaldson-Thomas bundle & Anti-Self-Dual
bundle\\\hline
\end{tabular}
\]

\begin{definition}
Suppose $M$ is a special Riemannian $\mathbb{A}$-manifold. A pair $\left(
C,D_{E}\right)  $ is called an (i) special $\mathbb{A}$-cycle if $C$ is an
$\mathbb{A}$-submanifold in $M$ and $D_{E}$ is a special $\mathbb{A}%
$-connection over $C$ or (ii) special $\frac{1}{2}\mathbb{A}$-Lagrangian cycle
(or type I or II) if $C$ is a special $\frac{1}{2}\mathbb{A}$-Lagrangian
submanifold (or type I or II) in $M$ and $D_{E}$ is a special $\frac{1}%
{2}\mathbb{A}$-connection over $C$.
\end{definition}

Notice that all special cycles are calibrated. For example a special
$\mathbb{C}$-cycle $\left(  C,D_{E}\right)  $ in a K\"{a}hler manifold is
calibrated by $\exp\omega$ because a complex submanifold $C$ of dimension $2k$
is volume calibrated by $\omega^{k}$ and a Hermitian Yang-Mills connection
over $C$ is Yang-Mills calibrated by $\omega^{k-2}$. The following table gives
the common names of special cycles and their calibrating forms.%

\[%
\begin{tabular}
[c]{|c||c||c|c|}\hline
$%
\begin{array}
[c]{l}%
\,\\
\,
\end{array}
$ & Special $\mathbb{A}$-cycle & Special $\frac{1}{2}\mathbb{A}$-Lagr. cycle,
I & Special $\frac{1}{2}\mathbb{A}$-Lagr. cycle, II\\
& (calibrating form) & (calibrating form) & (calibrating form)\\\hline\hline
&  Points or $M$ + &  & \\
$\mathbb{R}$ & Unitary flat bundle & n/a & n/a\\
$%
\begin{array}
[c]{l}%
\,\\
\,
\end{array}
$ & ($\exp\nu_{M}$) &  & \\\hline
&  Complex submanifold + & Lagr. submfd phase $0$ + & Lagr. submfd phase
$\pi/2$ +\\
$\mathbb{C}$ & Herm. Yang-Mills bundle & Unitary flat bundle & Unitary flat
bundle\\
$%
\begin{array}
[c]{l}%
\,\\
\,
\end{array}
$ & ($\exp\omega$) & ($\exp\left(  \operatorname{Re}\Omega\right)  $) &
($\exp\left(  \operatorname{Im}\Omega\right)  $)\\\hline
&  Quaternionic submfd + & Complex Lagr. submfd + & Complex Lagr. submfd +\\
$\mathbb{H}$ & ASD connection & Herm. Yang-Mills bdl & Herm. Yang-Mills bdl\\
$%
\begin{array}
[c]{l}%
\,\\
\,
\end{array}
$ & ($\exp\Theta$) & ($\exp\left(  \operatorname{Re}\Omega_{I}^{n}\right)  $
and $\exp\omega_{J}$) & ($\exp\left(  \operatorname{Re}\Omega_{I}^{n}\right)
$ and $\exp\omega_{J}$)\\\hline
&  The manifold $M$ + & Associative submfd + & Coassociative submfd +\\
$\mathbb{O}$ & $G_{2}$-Donaldson-Thomas bdl & Unitary flat bundle & ASD
connection\\
$%
\begin{array}
[c]{l}%
\,\\
\,
\end{array}
$ & ($\exp\Theta_{M}$) & ($\exp\Omega_{M}$) & ($\exp\Theta_{X}$)\\\hline
\end{tabular}
\]

(Special) $\mathbb{A}$-geometry studies (special) $\mathbb{A}$-cycles and
(special) $\frac{1}{2}\mathbb{A}$-Lagrangian geometry studies (special)
$\frac{1}{2}\mathbb{A}$-Lagrangian cycles on (special) Riemannian $\mathbb{A}$-manifolds.

Remark: (Special) $\mathbb{C}$-geometry and $\mathbb{R}$-Lagrangian geometry
are basically the complex algebraic geometry and the symplectic geometry.
Special $\mathbb{R}$-Lagrangian geometry is important in the SYZ mirror
conjecture. Special $\mathbb{C}$-Lagrangian geometry is studied in \cite{Le
Lag HK Legendre} and is closely related to the classical Plucker formula. For
(special) $\mathbb{O}$-manifolds, these geometries are discussed by Donaldson
and Thomas \cite{DT}, Hitchin, Gukov, Yau and Zaslow \cite{GYZ}, Lee and the
author \cite{LL}, \cite{Le Asym Cyl G2}.

\subsection{\label{1Sec Fouier A geom}Fourier transformation of $\mathbb{A}$-geometry}

First we review the Fourier transformation in geometry. Classically Fourier
transformation is a duality between functions on a vector space $V\cong
\mathbb{R}^{n}$ and on its dual vector space $V^{\ast}$. It is given by
\[
f\left(  x\right)  \rightarrow\hat{f}\left(  y\right)  =\frac{1}{\left(
2\pi\right)  ^{n}}\int_{V}f\left(  x\right)  e^{ix\cdot y}dx\text{.}%
\]

The Fourier transformation on the geometry of flat tori is similar. Suppose
$T=V/\Lambda$ is any n dimensional torus, i.e. $\Lambda\cong\mathbb{Z}^{n}$ is
a lattice in $V\cong\mathbb{R}^{n}$. The dual torus $T^{\ast}$ is defined as
$V^{\ast}/\Lambda^{\ast}$ where $\Lambda^{\ast}$ consists of those $\phi\in
V^{\ast}$ with $\phi\left(  \Lambda\right)  \subset\mathbb{Z}$. This
relationship is reflexive, $\left(  T^{\ast}\right)  ^{\ast}=T$. Moreover,
$T^{\ast}$ can be naturally identified with the moduli space of flat $U\left(
1\right)  $-connections on $T$,
\[
\mathcal{M}^{U\left(  1\right)  \text{-flat}}\left(  T\right)  \cong T^{\ast
}.
\]

On $T\times T^{\ast}$, there is \textit{universal Poincar\'{e} line bundle}
$\mathbf{L}$ with a universal connection, $\mathbf{D}=d+\pi i%
{\textstyle\sum_{j=1}^{n}}
\left(  y^{j}dy_{j}-y_{j}dy^{j}\right)  $, where $y_{j}$'s are coordinates on
$T^{\ast}$ dual to the linear coordinates $y^{j}$'s on $T$. Its curvature,
$\mathbf{F}=2\pi i\Sigma dy^{j}\wedge dy_{j}$ plays the role of the kernel
function in the classical Fourier transformation.

\bigskip

\textbf{Fourier transformation of }$\mathbb{R}$\textbf{-geometry on }$T$

On the topological level, the Fourier transformation is given by,
\begin{align*}
\mathcal{F}  &  :H^{k}\left(  T,\mathbb{Z}\right)  \overset{\simeq
}{\rightarrow}H^{n-k}\left(  T^{\ast},\mathbb{Z}\right) \\
\mathcal{F}\left(  \phi\right)   &  =\int_{T}\phi\wedge e^{\frac{i}{2\pi
}\mathbf{F}},
\end{align*}
and we also have a similar one for K-groups.

On the \textit{flat }level, we consider flat bundles over $T$ or points in
$T$, i.e. cycles $\left(  C,D_{E}\right)  $ that are calibrated by $\exp
\nu_{M}$, where $\nu_{M}$ is the volume form on $T$. Any flat $U\left(
r\right)  $-bundle $E$ over $T$ is isomorphic to an direct sum of flat line
bundles, unique up to permutations. This can be interpreted as an
identification of moduli spaces of special $\mathbb{R}$-cycles on $T$ and on
$T^{\ast}$, via the Fourier transformation.

\bigskip

\textbf{Fourier transform of }$\mathbb{C}$\textbf{-geometry on }$T$

When $T$ and $T^{\ast}$ are Abelian varieties, Mukai shows that the Fourier
transformation,
\[
\mathcal{F}\left(  \cdot\right)  =\mathbf{R}\pi_{1\ast}\left(  \pi_{2}^{\ast
}\left(  \cdot\right)  \otimes\mathcal{L}\right)
\]
is an equivalence of derived categories of coherent sheaves, $D^{b}\left(
T\right)  $ and $D^{b}\left(  T^{\ast}\right)  $, with the inversion property.
This has far reaching implications in the theory of Abelian varieties.

\bigskip

\textbf{Fourier transform of }$\mathbb{H}$\textbf{-geometry on }$T$

Unforturnately, we only know of a low dimension example in this case: A flat
torus $T$ of dimension four is a special $\mathbb{H}$-manifold. Braam and
Schenk (see e.g. \cite{DT}) show that the Fourier transformation of an ASD
connection over $T$ without any flat factor is another ASD connection over
$T^{\ast}$. Moreover this bijection between their moduli spaces is an isometry
with respect to the Weil-Peterrson L$^{2}$-metrics.

\bigskip

In summary we expect that the Fourier transform on flat $\mathbb{A}$-tori
gives a correspondence:%

\[
\mathbb{A}\text{-Geometry}\left(  T\right)  \longleftrightarrow\mathbb{A}%
\text{-Geometry}\left(  T^{\ast}\right)  .
\]

\subsection{\label{1Sec B Lag fibratn}$\frac{1}{2}\mathbb{A}$-Lagrangian fibrations}

On a Riemannian $\mathbb{A}$-manifold $M$, a $\frac{1}{2}\mathbb{A}%
$-\textit{Lagrangian fibrations}
\[
f:M\rightarrow B
\]
is a smooth map $f$ such that its generic fibers are smooth $\frac{1}%
{2}\mathbb{A}$-Lagrangian submanifolds of $M$. We usually also assume that $f$
has a section which is also a $\frac{1}{2}\mathbb{A}$-Lagrangian submanifold.

$\frac{1}{2}\mathbb{A}$-Lagrangian fibrations have been playing important
roles in various branches of mathematics and physics. In symplectic geometry,
a Lagrangian fibration with a Lagrangian section is a very important
structure. It is sometimes called a completely integrable system, or a real
polarization when we try to quantize the symplectic manifold. Familiar
examples include toric varieties.

In string theory, Strominger, Yau and Zaslow propose that mirror symmetry
should be explained in terms of the fiberwise Fourier transformation along
special Lagrangian fibrations on mirror Calabi-Yau manifolds \cite{SYZ}.

Heuristically, one should view a $\frac{1}{2}\mathbb{A}$-Lagrangian fibration
with a section on $M$ as a \textit{global de-complexification} of $M$.

\bigskip

Remark: When $\mathbb{A}$ equals $\mathbb{C}$ or $\mathbb{H}$, a generic fiber
of any $\frac{1}{2}\mathbb{A}$-Lagrangian fibration is a torus. However in a
$Spin\left(  7\right)  $-manifold or a $G_{2}$-manifold, such a fiber is
expected to be either a torus or a K3 surface, at least near a adiabatic limit
(see \cite{LL}).

\subsection{\label{1Sec Mirror dual}Mirror duality between different geometries}

Following SYZ proposal, we would be interested in a fiberwise Fourier
transformation on any Riemannian $\mathbb{A}$-manifold $M$ with a given
$\frac{1}{2}\mathbb{A}$-Lagrangian fibration and a section. First we want to
construct a dual torus fibration
\[
g:W\rightarrow B.
\]
First this would require each fiber of $f$ to have the same dimension.
Typically, $M$ would then be a special $\mathbb{A}$-manifold, except in the
octonionic case. Second, in order the perform the Fourier transformation on
fibers, we require these fiber tori to be flat, at least in the limit. To be
more precise, we assume that $M$ has an one parameter family of such metrics
parametrized by $t\in\lbrack0,\infty)$, such that as $t\rightarrow\infty$ the
second fundamental form of each smooth fiber goes to zero. This is related to
the so-called \textit{large structure limit }in the physics literature. Then
one expects that (i) the total space of the dual torus fibration is also a
special $\mathbb{A}$-manifold and (ii) the fiberwise Fourier transformation
will give an equivalence of geometries:
\[
\text{(Special) }\mathbb{A}\text{-Geometry}\left(  M\right)
\longleftrightarrow\text{(Special) }\frac{1}{2}\mathbb{A}\text{-Lagrangian
Geometry}\left(  W\right)  .
\]
Moreover the relationship between $M$ and $W$ should be reflexive.

The above picture has been over-simplified. There are many subtleties
involved. Many of them are related to quantum corrections, an issue we have
not addressed here. The most famous example is the mirror symmetry conjecture
for special $\mathbb{C}$-manifolds (i.e. Calabi-Yau manifolds). The fiberwise
Fourier transformation in this case has been studied by many people including
Gross, Hitchin, Kontsevich, Ruan, Vafa, Witten, Yau, Zaslow, the author and
many others.

We would indicate how the fiberwise Fourier transformation works in the
simplest situation and show how the symplectic structure (determined by
$\exp\omega$) is being transformed to the complex structure (determined by
$\Omega$) on the mirror: Suppose $M=\mathbb{C}^{n}=\mathbb{R}^{n}%
\times\mathbb{R}^{n}$ with the standard complex structure $J$, i.e.
$\Omega_{M}=dz^{1}\wedge\cdots\wedge dz^{n}$, and with a symplectic structure
$\omega_{M}=\Sigma\phi_{ij}\left(  x\right)  dx^{i}dy^{j}$ which is invariant
under translations along $y$ directions. Then $W=\mathbb{R}^{n}\times
\mathbb{R}^{n\ast}$ with coordinates $x^{i}$'s and $y_{i}$'s. By direct
computations, under the fiberwise Fourier transformation on differential
forms,
\[
\left(  \cdot\right)  \rightarrow\int\left(  \cdot\right)  e^{\mathbf{F}%
}dy^{1}dy^{2}\cdots dy^{n}.
\]
where $\mathbf{F}=i\Sigma dy^{j}\wedge dy_{j}$ is the universal curvature
form, we have
\begin{align*}
\exp\left(  \omega_{M}\right)   &  \rightarrow\Omega_{W}=\prod\left(
\phi_{ij}dx^{i}+idy_{j}\right) \\
\Omega_{M}  &  \rightarrow\exp\left(  \omega_{W}\right)  =\exp\left(
{\textstyle\sum}
dx^{i}\wedge dy_{i}\right)  \text{.}%
\end{align*}
Thus we see how variation of symplectic structures on $M$ corresponds to
variation of complex structures on $W$ explicitly (see \cite{Le MS corr} for
more details).

\section{\label{1Sec Rem Qu}Remarks and questions}

In this last section, we remark on some other aspects of geometry over
$\mathbb{A}$ with $\mathbb{A}=\mathbb{R}$, $\mathbb{C}$, $\mathbb{H}$,
$\mathbb{O}$, and mention a few interesting questions.

\bigskip

\textbf{Triality transformation}

As we discussed in section \ref{1Sec Mirror dual}, mirror symmetry for
Calabi-Yau manifolds and hyperk\"{a}hler manifolds is a duality transformation
for the geometry of these manifolds. In its simplest form, it can be viewed as
the duality between a vector space and its dual vector space. The novelty
about octonion in algebra is the \textit{triality} (see e.g. \cite{Baez}). It
would be important to introduce the triality transformation to the geometry
for $\mathbb{O}$-manifolds. The work of Gukov, Yau and Zaslow \cite{GYZ} might
be related to this.

\bigskip

\textbf{Moduli space of connections}

In gauge theory, we consider the space of connections $\mathcal{A}$ on a
bundle $E$ modulo the group of gauge transformations $\mathcal{G}=Aut\left(
E\right)  $, and call this the moduli space of connections on $E$. On a
K\"{a}hler manifold $M$, the moduli space of special $\mathbb{C}$-connections
(i.e. Hermitian Yang-Mills connections) is the symplectic reduction $\mu
^{-1}\left(  0\right)  /\mathcal{G}$ of the space of $\mathbb{C}$-connections
by $\mathcal{G}$, because
\begin{align*}
\mu &  :\mathcal{A}\rightarrow\Omega^{2n}\left(  M,ad\left(  E\right)  \right)
\\
\mu\left(  D_{E}\right)   &  =F_{E}\wedge\omega^{n-1}%
\end{align*}
is the moment map for the action of $\mathcal{G}$ on $\mathcal{A}$. By the
work of Donaldson, Uhlenbeck and Yau, we can also view this as the complex
quotient of the space of $\mathbb{C}$-connections by $\mathcal{G}^{\mathbb{C}%
}$. Similarly on a hyperk\"{a}hler manifold $M$, there is a hyperk\"{a}hler
moment map,
\begin{align*}
\mu &  :\mathcal{A}\rightarrow\Omega^{4n}\left(  M,ad\left(  E\right)
\right)  \otimes\mathbb{R}^{3}\\
\mu\left(  D_{E}\right)   &  =\left(  F_{E}\wedge\omega_{I}^{2n-1},F_{E}%
\wedge\omega_{J}^{2n-1},F_{E}\wedge\omega_{K}^{2n-1}\right)  ,
\end{align*}
such that the moduli space of special $\mathbb{H}$-connections on $E$ is the
hyperk\"{a}hler quotient of $\mathcal{A}$ by $\mathcal{G}$, restricted to the
space of $\mathbb{H}$-connections. On a $G_{2}$-manifold $M=X\times S^{1}$,
one does not have the notion of an octonionic quotient. However we still can
define $\mu$ in the same way and its zeros correspond to connections pullback
from $X$. Therefore we have a similar picture as before, namely the space of
$Spin\left(  7\right)  $-Donaldson-Thomas connections on $M$ with $\mu=0$
quotienting by $\mathcal{G}$ is the moduli space of $G_{2}$-Donaldson-Thomas
connections on $X$.

$\bigskip$

$\mathbb{A}$\textbf{-manifolds without metrics}

Real and complex manifolds are usually defined using coordinate charts and
requiring their transition functions to be differentiable and holomorphic
respectively. However such definitions can not be generalized to the
quaternionic and octonionic cases. For example any smooth map which preserves
quaternionic structures must be affine (e.g. \cite{Besse}). The correct
generalization is to use a torsion free connection on $M$ preserving an
$\mathbb{A}$-structure on the frame bundle: We denote by $\mathcal{G}%
_{\mathbb{A}}\left(  n\right)  $ the group of twisted isomorphisms $\phi$ of
$\mathbb{A}^{n}$, but do not require $\phi$ to be an isometry. Similarly, we
have $\mathcal{H}_{\mathbb{A}}\left(  n\right)  $ for special twisted isomorphisms.

\begin{definition}
A smooth manifold $M$ is called a (special) $\mathbb{A}$-manifold if the
structure group of its frame bundle has a reduction to $\mathcal{G}%
_{\mathbb{A}}\left(  n\right)  $ (resp. $\mathcal{H}_{\mathbb{A}}\left(
n\right)  $) together with a torsion free connection.
\end{definition}

The following table gives explicit descriptions of $\mathcal{G}_{\mathbb{A}%
}\left(  n\right)  $ and $\mathcal{H}_{\mathbb{A}}\left(  n\right)  $,
together with the usual names for these manifolds.%

\[%
\begin{tabular}
[c]{|l|l|l|}\hline
& $\mathcal{G}_{\mathbb{A}}\left(  n\right)
\begin{array}
[c]{l}%
\,\\
\,
\end{array}
$ & $\mathcal{H}_{\mathbb{A}}\left(  n\right)  $\\
& ($\mathbb{A}$-manifolds) & (Special $\mathbb{A}$-manifolds)\\\hline\hline
$\mathbb{R}$ & $GL\left(  n,\mathbb{R}\right)
\begin{array}
[c]{l}%
\,\\
\,
\end{array}
$ & $GL^{+}\left(  n,\mathbb{R}\right)  $\\
& (Manifolds) & (O$\text{riented manifolds}$)\\\hline\cline{3-3}%
$\mathbb{C}$ & $GL\left(  n,\mathbb{C}\right)
\begin{array}
[c]{l}%
\,\\
\,
\end{array}
$ & $SL\left(  n,\mathbb{C}\right)  $\\
& (C$\text{omplex manifolds}$)\thinspace & \\\hline\cline{3-3}%
$\mathbb{H}$ & $GL\left(  n,\mathbb{H}\right)  \mathbb{H}^{\times}
\begin{array}
[c]{l}%
\,\\
\,
\end{array}
$ & $GL\left(  n,\mathbb{H}\right)  $\\
& (Q$\text{uaternionic manifolds}$) & (H$\text{ypercomplex manifolds}%
$)\\\hline\cline{3-3}%
$\mathbb{O}$ & $Spin\left(  7\right)
\begin{array}
[c]{l}%
\,\\
\,
\end{array}
$ & $G_{2}$\\
& ($Spin\left(  7\right)  \text{-manifolds}$) & ($G_{2}\text{-manifolds}%
$)\\\hline\cline{3-3}%
\end{tabular}
\]

\bigskip

$\mathbb{A}$\textbf{-torsions}

Ray and Singer define analytic torsions for real and complex manifolds, they
are important non-local invariants for these manifolds and play important
roles in the family index theory. The quaternionic generalization of the
analytic torsion is discussed by Leung, Yi \cite{Le Yi torsion} and Koehler,
Weingartmath \cite{KM Quat torsion}. It is natural to ask for their analog in
the octonionic case and their common roles in the geometry over $\mathbb{A}$.

\bigskip

$G_{2}$\textbf{-symplectic manifolds}

A different generalization of Riemannian $\mathbb{A}$-manifolds is by not
requiring the almost complex structures on $M$ to be integrable, but only
require the \textit{closedness} of the calibrating differential forms like the
K\"{a}hler form or the holomorphic volume form. It turns out that
integrability is automatic in the Calabi-Yau case, hyperk\"{a}hler case and
$Spin\left(  7\right)  $-case. And we only obtain two new classes of manifolds
this way, namely (1) symplectic manifolds and (2) almost $G_{2}$-manifolds. In
\cite{Le TQFT} the author discuss TQFT on almost $G_{2}$-manifolds.

Notice that $\frac{1}{2}\mathbb{A}$-Lagrangian submanifolds can be defined
using vanishing of certain two forms, thus they are well-defined for such
manifolds. In the symplectic case, they are simply Lagrangian submanifolds in
the usual sense. Indeed, if $M$ is an almost $G_{2}$-manifold, i.e. there
exists a closed non-degenerate three form on $M$, then the space
$\mathcal{L}M$ of all unparametrized loop in $M$ has a natural symplectic
structure, as observed by Movshev.

\bigskip

\textbf{Twistor theory for octonionic manifolds}

On any oriented Riemannian four manifold $M$, we can define a\textit{\ twistor
space} $Z$ with a fibration
\[
f:Z\rightarrow M
\]
whose fiber over $x$ is the $S^{2}$-family of Hermitian complex structures on
$T_{x}M$. The twistor space $Z$ has a natural almost complex structure which
is integrable if and only if the Weyl curvature tensor of $M$ is self-dual.
Penrose's twistor transformation gives a correspondence between the conformal
geometry of $M$ and the complex geometry of $Z$ (see e.g. \cite{Besse}). There
is a natural generalization of the twistor transform for any quaternionic
K\"{a}hler manifold, which is a transformation between the quaternionic
geometry of $M$ and the complex geometry of its twistor space $Z$.

In certain sense, the twistor theory is another form of a
\textit{de-complexification} of $M$. Therefore it is natural to ask whether
there is an analog theory for $G_{2}$-manifolds, or even $Spin\left(
7\right)  $-manifolds. One would expect that such twistor spaces would be
quaternionic manifolds.

\bigskip

\textit{Acknowledgments: This paper is partially supported by NSF/DMS-0103355.
The author expresses his gratitude to S.L. Kong, J.H. Lee for useful discussions.}

\bigskip\

\bigskip

Address: School of Mathematics, University of Minnesota, Minneapolis, MN
55454, USA.

Email: LEUNG@MATH.UMN.EDU

\bigskip

\bigskip
\end{document}